\input amstex


\def\b1{\text{\bf 1}}

\def\BC{{\Bbb C}}

\def\BQ{{\Bbb Q}}

\def\btau{\bar{\tau}}

\def\BZ{{\Bbb Z}}
\def\CA{{\Cal A}}
\def\CB{{\Cal B}}

\def\CD{{\Cal D}}

\def\CO{{\Cal O}}

\def\CT{{\Cal T}}

\def\dpar{\partial}
\def\dplus{\buildrel\cdot\over{+}}

\def\fA{{\frak A}}
\def\fa{{\frak a}}
\def\fB{{\frak B}}
\def\fb{{\frak b}}
\def\fC{{\frak C}}
\def\fD{{\frak D}}
\def\fE{{\frak E}}
\def\fg{{\frak g}}


\def\btu{\bigtriangleup}
\def\hra{\hookrightarrow}
\def\iso{\buildrel\sim\over\longrightarrow} 

\def\lra{\longrightarrow}

\parskip=6pt

\documentstyle{amsppt}
\document
\NoBlackBoxes


\centerline{\bf Gerbes of chiral differential operators. III}

\bigskip
\centerline{Vassily Gorbounov, Fyodor Malikov, Vadim Schechtman}
\bigskip


\bigskip\bigskip 

\centerline{\bf Introduction}

\bigskip\bigskip

This note is a sequel to [GII]. 
Its aim is to "switch on an exterior vector bundle" 
into the framework of {\it op. cit.} 

Let $X$ be a smooth scheme over a fixed ground 
ring $k$ containing $1/2$ and $E$ be a vector bundle 
(i.e. a locally free $\CO_X$-module) of finite rank over $X$. 
Consider the exterior algebra $\Lambda E=\oplus_{i=0}^{rk(E)}\ 
\Lambda^i E$ (over $\CO_X$); this is a sheaf of commutative superalgebras 
over $X$,  
where by definition $\CO_X$ is purely even and the parity of a component 
$\Lambda^i E$ is equal to the parity of $i$. 

In this note we study the chiral counterparts of the sheaf $\CD_{\Lambda E}$ 
of superalgebras of differential operators acting on $\Lambda E$. 
Similarly to {\it op. cit.}, these {\it chiral sheaves of differential 
operators on 
$\Lambda E$} exist locally and are by no means unique; the corresponding 
categories form a {\it champ en groupoids} $\fD_{\Lambda E}$ over $X$, 
called the {\it gerbe of chiral differential operators on $\Lambda E$}. 

Our first main result (see Theorem 5.9) says that the characteristic 
class $c(\fD_{\Lambda E})$ lies in the second hypercohomology group 
$H^2(X;\Omega_X^{[2,3\rangle})$ (i.e. in the same group where  
$c(\fD_X)$ lies) and is equal to
$$
c(\fD_{\Lambda E})=c(\Theta_{X/k})-c(E)=c(\Omega^1_{X/k})-c(E)
\eqno{(0.1)}
$$
where $c(E)$ is the "Atiyah-Chern-Simons" class defined in [GII], 7.6. 
Here $\Theta_{X/k}$ is the tangent bundle. Recall that 
$\Omega^{[2,3\rangle}_X$ denotes the length $1$ complex of sheaves 
$\Omega^2_{X/k}\lra\Omega^{3,closed}_{X/k}$, with $\Omega^2_{X/k}$ living 
in degree $0$. 

As usually, we obtain in fact a stronger statement, namely the equality 
(0.1) "on the level of cocycles". As a corollary of this, we 
conclude that for $E=\Theta_{X/k}$ or $E=\Omega^1_{X/k}$ our 
gerbes admit a canonical global section. In other words, 
there exist canonically defined {\it the} sheaves of chiral do 
$\CD^{ch}_{\Lambda \Theta_X}$ and $\CD^{ch}_{\Omega^{\cdot}_X}$. 

Section 6 is devoted to the study of the last sheaf, 
which is nothing but (the underlying sheaf of)  
{\it chiral de Rham complex} from [MSV]. We obtain the transformation 
laws of $4$ local generators of $N=2$ supersymmetry $Q, J, G$ and $L$, 
see Theorem 6.25. In particular, the component $Q_0$ of the field 
$Q(z)$ is a globally defined square zero derivation of 
$\CD^{ch}_{\Omega^{\cdot}_X}$, which is the {\it chiral de Rham 
differential} from {\it op. cit.} 

This completes an alternative construction of the chiral de Rham 
complex sketched in Section 6 of [MSV]. Its difference from the 
original construction is that it does not use Wick theorem and the arguments 
of "formal geometry". 

In the last section we show that as a simple consequence of 
the Poincar\'e-Birkhoff-Witt theorem for $\CD^{ch}_{\Omega^{\cdot}_X}$ 
and the Lefschetz fixed point theorem one gets a "moonshine style" 
formula, cf. Theorem 7.9.  

The work was done while V.S. visited IHES. He is grateful to the Institute 
for the excellent working atmosphere.

\bigskip\bigskip


\centerline{\bf \S 1. Preliminaries}

\bigskip\bigskip

{\bf 1.1.} We keep the assumptions of [GII]. We assume that the ground ring $k$ 
contains $1/2$. For a $k$-supermodule $M$, we denote by $M^{ev}$ 
(resp. $M^{odd}$) the submodule of even (resp. odd) elements, 
so that $M=M^{ev}\oplus M^{odd}$. For a homogeneous element $a\in M$, 
we denote by $p(a)\in\BZ/2\BZ$ its parity.  
When we speak about graded $k$-supermodules  
$M=\oplus_{i\in I}\ M_i$ we imply that the $I$-grading is compatible 
with the parity, i.e. $M^x=\oplus_i\ M_i^x$ where $M_i^x=M^x\cap M_i,\ 
x=ev$ or $odd$. 

Let $A$ be a commutative $k$-superalgebra. A {\it Lie superalgebroid} over $A$  
is a Lie superalgebra over $k$ equipped with a structure of an $A$-module, 
such that the identities [GII] (0.2.1) and (0.2.2) hold true. 

{\bf 1.2.} A {\it $\BZ_{\geq 0}$-graded vertex superalgebra} (over $k$) is a 
$\BZ_{\geq 0}$-graded $k$-supermodule $V=\oplus_{i\geq 0}\ V_i$  equipped 
with a distinguished even vector $\b1\in V_0$ ({\it vacuum vector}) 
and a family of bilinear operations 
$$
_{(n)}:\ V\times V\lra V,\ n\in\BZ,
$$
such that 
$$
p(a_{(n)}b)=p(a)+p(b);\ \ V_{i(n)}V_j\subset V_{i+j-n-1}
\eqno{(1.2.1)}
$$
The following properties must hold:
$$
\b1_{(n)}a=\delta_{n,-1}a;\ a_{(n)}\b1=0\ \text{for }n\geq 0,\ 
a_{(-1)}\b1=a
\eqno{(1.2.2)}
$$
and 
$$
\sum_{j=0}^\infty\ \binom{m}{j}(a_{(n+j)}b)_{(m+l-j)}c=
$$
$$
=
\sum_{j=0}^\infty\ (-1)^j\binom{n}{j} 
\bigl\{a_{(m+n-j)}b_{(l+j)}c-(-1)^{n+p(a)p(b)}b_{(n+l-j)}a_{(m+j)}c\bigr\}
\eqno{(1.2.3)}
$$
for all $m,n,l\in \BZ$. A particular case of 
(1.2.3) corresponding to $m=0$: 
$$
(a_{(n)}b)_{(l)}c=\sum_{j=0}^\infty\ (-1)^j\binom{n}{j}
\bigl\{ a_{(n-j)}b_{(l+j)}c-(-1)^{n+p(a)p(b)}b_{(n+l-j)}a_{(j)}c\bigr\}
\eqno{(1.2.4)}
$$
Setting $n=l=-1$ we get 
$$
(a_{(-1)}b)_{(-1)}c=\sum_{j=0}^\infty\ 
\bigl\{ a_{(-1-j)}b_{(-1+j)}c+(-1)^{p(a)p(b)}b_{(-2-j)}a_{(j)}c\bigr\}
\eqno{(1.2.5)}
$$

In the sequel we shall work only with $\BZ_{\geq 0}$-graded vertex 
superalgebras, and call them simply vertex superalgebras. 
This $\BZ_{\geq 0}$-grading will be called the grading by {\it 
conformal weight}. 

{\bf 1.3.} Let $V$ be a vertex superalgebra.   
The even operators $\dpar^{(j)}:\ V\lra V$ of degree $j$ ($j\in\BZ_{\geq 0}$) 
are defined in the same manner as in [GII], (0.5.5), and they satisfy 
[GII], (0.5.7), (0.5.8), (0.5.10) and (0.5.11). The "supercommutativity" 
formula reads as 
$$
a_{(n)}b=(-1)^{n+p(a)p(b)+1}\sum_{j\geq 0}\ (-1)^j \dpar^{(j)}(b_{(n+j)}a)
\eqno{(1.3.1)}
$$
and we have the usual OPE formula
$$
[a_{(m)},b_{(n)}]=\sum_{j\geq 0}\ \binom{m}{j}(a_{(j)}b)_{(m+n-j)}
\eqno{(1.3.2)}
$$
where in the left hand side stands the supercommutator 
$$
[a_{(m)},b_{(n)}]:=a_{(m)}b_{(n)}-(-1)^{p(a)p(b)}b_{(n)}a_{(m)}
\eqno{(1.3.3)}
$$ 

\bigskip\bigskip

\centerline{\bf \S 2. Vertex Superalgebroids} 

\bigskip\bigskip

{\bf 2.1.} An {\it extended Lie superalgebroid} (over $A$) is a quintuple 
$\CT=(A,T,\Omega,\dpar,\langle,\rangle)$ where $T$ is a Lie superalgebroid 
over $A$, $\Omega$ is an $A$-module equipped with a structure of 
a module over the Lie superalgebra $T$, $\dpar:\ A\lra\Omega$ an even 
$A$-derivation and a morphism of $T$-modules, 
$\langle,\rangle:\ T\times\Omega\lra A$ is an even $A$-bilinear pairing. 

The following identities must be true $(a\in A, \tau,\nu\in T, 
\omega\in\Omega)$:

$$
\langle\tau,\dpar a\rangle=\tau(a)
\eqno{(2.1.1)}
$$ 
$$
\tau(a\omega)=\tau(a)\omega+(-1)^{p(\tau)p(a)}a\tau(\omega)
\eqno{(2.1.2)}
$$
$$
(a\tau)(\omega)=a\tau(\omega)+\langle\tau,\omega\rangle\dpar a
\eqno{(2.1.3)}
$$
$$
\tau(\langle\nu,\omega\rangle)=\langle [\tau,\nu],\omega\rangle+ 
(-1)^{p(\tau)p(\nu)}\langle\nu,\tau(\omega)\rangle
\eqno{(2.1.4)}
$$

For example to a Lie superalgebroid $T$ one associates canonically 
an extended superalgebroid with $\Omega=Hom_A(T,A)$, as in 
[GII], 1.2.     

{\bf 2.2. De Rham - Chevalley complex.} Let $\CT=(A,T,\Omega,\ldots)$ be 
an extended Lie $A$-superalgebroid. Let us define $A$-modules $\Omega^i=
\Omega^i(\CT),\ 
i\in\BZ_{\geq 0}$, as follows. Set $\Omega^0=A,\ \Omega^1=\Omega$. For 
$i\geq 2$, $\Omega^i$ is the submodule of the module of  
$A$-polylinear 
homomorphisms $h$ from $T^{i-1}$ to $\Omega$ such that the function 
$\langle\tau_1,h(\tau_2,\ldots,\tau_i)\rangle$ is skew symmetric 
(in the graded sense) with 
respect to all permutations of $(\tau_1,\ldots,\tau_i)$. 

For example, if $\CT$ is associated with a Lie superalgebroid $T$ as above 
then $\Omega^i=Hom_A(\Lambda^i_AT,A)$.  

Let us define the even maps $d_{DR}=d_{DR}^i:\ \Omega^i\lra\Omega^{i+1}$ 
as follows. 
For $i=0$ we set $d_{DR}a=-\dpar a$. For $i\geq 1$ we set  
$$
d_{DR}h(\tau_1,\ldots,\tau_i)=d_{Lie}h(\tau_1,\ldots,\tau_i)-
(-1)^{p(h)p(\tau_1)}\dpar\langle\tau_1,h(\tau_2,\ldots,\tau_i)\rangle
\eqno{(2.2.1)}
$$
where
$$
d_{Lie}h(\tau_1,\ldots,\tau_i)=\sum_{j=1}^i\ (-1)^{j+1+p(\tau_j) 
(p(\tau_1)+\ldots+p(\tau_{j-1}))}
\tau_j(h(\tau_1,\ldots,\widehat{\tau}_j,\ldots,\tau_i))+
$$
$$
+\sum_{1\leq j<l\leq i}\ (-1)^{j+l+p(\tau_j)(p(\tau_1)+\ldots+p(\tau_{j-1}))+
p(\tau_l)(p(\tau_1)+\ldots+\widehat{p(\tau_j)}+\ldots+p(\tau_{l-1}))}\times 
$$
$$ 
\times h([\tau_j,\tau_l],\tau_1,\ldots,\widehat{\tau}_j,\ldots,
\widehat{\tau}_l,\ldots,\tau_i)
\eqno{(2.2.2)}
$$
For example, 
$$
d_{DR}\omega(\tau)=(-1)^{p(\omega)p(\tau)}\bigl\{\tau(\omega)-
\dpar\langle\tau,\omega\rangle\bigr\},\ 
\eqno{(2.2.3)}
$$
for $\omega\in\Omega^1=\Omega$; and 
$$
d_{DR}h(\tau_1,\tau_2)=-h([\tau_1,\tau_2])+
(-1)^{p(h)p(\tau_1)}\tau_1(h(\tau_2))- 
$$
$$ 
-(-1)^{p(\tau_2)(p(h)+p(\tau_1)}\tau_2(h(\tau_1))-
(-1)^{p(h)p(\tau_1)}\dpar\langle\tau_1,h(\tau_2)\rangle,\ 
\eqno{(2.2.4)}
$$
for $h\in\Omega^2$. 

Let us introduce the action of the Lie algebra $T$ on the modules 
$\Omega^i$ by 
$$
\tau(h)(\tau_1,\ldots,\tau_{i-1})=\tau(h(\tau_1,\ldots,\tau_{i-1}))- 
$$
$$
-\sum_{j=1}^{i-1}\ (-1)^{p(\tau)(p(\tau_1)+\ldots+p(\tau_{j-1}))}
h(\tau_1,\ldots,[\tau,\tau_j],\ldots,\tau_i)
\eqno{(2.2.5)}
$$
Let us define the convolution operators $\langle\tau,\cdot\rangle:\ 
\Omega^i\lra\Omega^{i-1}$ by
$$
\langle\tau,h\rangle(\tau_1,\ldots,\tau_{i-2}\rangle=
(-1)^{p(\tau)p(h)}h(\tau,\tau_1,\ldots,\tau_{i-2})
\eqno{(2.2.6)}
$$
The maps $\{d_{DR}^i\}$ may be characterized as a unique collection 
of maps such that $d_{DR}^0=-\dpar$ and the {\it Cartan formula} 
$$
\tau(h)=\langle\tau,d_{DR}h\rangle+d_{DR}\langle\tau,h\rangle
\eqno{(2.2.7)}
$$
holds true. 

The maps $d_{DR}$ commute with the action of $T$. One checks that 
$d_{DR}^2=0$, so we get a complex $(\Omega^\bullet(\CT),d_{DR})$ 
called the {\it de Rham-Chevalley complex} of $\CT$.

{\bf 2.3.} A {\it vertex superalgebroid} is a septuple $\CA=(A,T,\Omega,\dpar,
\gamma,\langle\ ,\ \rangle,c)$ where $A$ is a supercommutative $k$-algebra, 
$T$ is a Lie superalgebroid over $A$, $\Omega$ is an $A$-module equipped 
with an action of the Lie superalgebra $T$, $\dpar:\ A\lra\Omega$ is an even  
derivation commuting with the $T$-action, 
$$
\langle\ ,\ \rangle:\ (T\oplus\Omega)\times (T\oplus\Omega)\lra A
$$
is a supersymmetric even $k$-bilinear pairing equal to zero on 
$\Omega\times\Omega$ and such 
that $\CT_{\CA}=(A,T,\Omega,\dpar,\langle\ ,\ \rangle|_{T\times\Omega})$ 
is an extended Lie superalgebroid over $A$; $c:\ T\times T\lra\Omega$ is a 
skew supersymmetric even $k$-bilinear pairing and 
$\gamma:\ A\times T\lra\Omega$ is an even $k$-bilinear map. 

The following axioms must hold $(a,b\in A;\ \tau,\tau_i\in T)$: 
$$
\gamma(a,b\tau)=\gamma(ab,\tau)-a\gamma(b,\tau)-
$$
$$
-(-1)^{p(\tau)(p(a)+p(b))}\tau(a)\dpar b-
(-1)^{p(a)p(b)+p(\tau)p(a)+p(\tau)p(b)}\tau(b)\dpar a
\eqno{(A1)}
$$
$$
\langle a\tau_1,\tau_2\rangle=a\langle\tau_1,\tau_2\rangle+
\langle\gamma(a,\tau_1),\tau_2\rangle-
(-1)^{p(a)(p(\tau_1)+p(\tau_2))}\tau_1\tau_2(a)
\eqno{(A2)}
$$
$$
c(a\tau_1,\tau_2)=ac(\tau_1,\tau_2)+\gamma(a,[\tau_1,\tau_2])-
$$
$$
-(-1)^{p(\tau_2)(p(\tau_1)+p(a))}\gamma(\tau_2(a),\tau_1)+
(-1)^{p(\tau_2)(p(\tau_1)+p(a))}\tau_2(\gamma(a,\tau_1))-
$$
$$
-(-1)^{p(a)(p(\tau_1)+p(\tau_2))}\frac{1}{2}\langle\tau_1,\tau_2\rangle\dpar a+
(-1)^{p(a)(p(\tau_1)+p(\tau_2))}\frac{1}{2}\dpar\tau_1\tau_2(a)-
$$
$$
-(-1)^{p(\tau_2)(p(a)+p(\tau_1))}\frac{1}{2}
\dpar\langle\tau_2,\gamma(a,\tau_1)\rangle
\eqno{(A3)}
$$
$$
\langle [\tau_1,\tau_2],\tau_3\rangle+
(-1)^{p(\tau_1)p(\tau_2)}\langle\tau_2,[\tau_1,\tau_3]\rangle=
\tau_1(\langle\tau_2,\tau_3\rangle)-
$$
$$
-(-1)^{p(\tau_1)p(\tau_2)}\frac{1}{2}\tau_2(\langle\tau_1,\tau_3\rangle)-
(-1)^{p(\tau_3)(p(\tau_1)+p(\tau_2))}
\frac{1}{2}\tau_3(\langle\tau_1,\tau_2\rangle)+
$$
$$
+(-1)^{p(\tau_1)p(\tau_2)}\langle\tau_2,c(\tau_1,\tau_3)\rangle+
(-1)^{p(\tau_3)(p(\tau_1)+p(\tau_2))}\langle\tau_3,c(\tau_1,\tau_2)\rangle
\eqno{(A4)}
$$
$$
d_{Lie}c(\tau_1,\tau_2,\tau_3)=-\frac{1}{2}\dpar\bigl\{
\langle [\tau_1,\tau_2],\tau_3\rangle+
(-1)^{p(\tau_2)p(\tau_3)}\langle [\tau_1,\tau_3],\tau_2\rangle-
$$
$$
-(-1)^{p(\tau_1)(p(\tau_2)+p(\tau_3))}\langle [\tau_2,\tau_3],\tau_1\rangle
-\tau_1(\langle\tau_2,\tau_3\rangle)+
(-1)^{p(\tau_1)p(\tau_2)}\tau_2(\langle\tau_1,\tau_3\rangle)-
$$
$$
-(-1)^{p(\tau_3)(p(\tau_1)+p(\tau_2))}2\langle\tau_3,c(\tau_1,\tau_2)\rangle
\bigr\}
\eqno{(A5)}
$$
where $d_{Lie}$ is defined by (2.2.2). 

{\bf 2.4.} All the constructions of [GII] generalize to the $\BZ/(2)$-graded 
case in an obvious manner.   

\bigskip\bigskip

\centerline{\bf \S 3. Some formulas}

\bigskip\bigskip

{\bf 3.1.} Let $A$ be a smooth $k$-algebra of relative dimension $n$, 
such that the $A$-module 
$T=Der_k(A)$ is free and admits a base $\{\btau_i\}$ consisting of commuting 
vector fields. Let $E$ be a free $A$-module of rank $m$, with a base 
$\{\phi_{\alpha}\}$. 
We shall call the set $\fg=\{\btau_i;\ \phi_{\alpha}\}\subset A\oplus E$ 
a {\it frame} of $(A,E)$. 

Consider a commutative $A$-superalgebra 
$\Lambda E=\oplus_{i=0}^m\ \Lambda_A^i(E)$ where the parity of 
$\Lambda^i_A(E)$ is equal to the parity of $i$. Each frame $\fg$ as above 
gives rise  to a $\Lambda E$-base $\{\tau_i;\ \psi_{\alpha}\}$  
of the Lie superalgebroid $T_{\Lambda E}=Der_k(\Lambda E)$, defined as follows. 
We extend the fields $\btau_i$ to derivations $\tau_i$ of the whole superalgebra 
$\Lambda E$ by the rule 
$$
\tau_i(a)=\btau_i(a);\ 
\tau_i(\sum\ a_{\alpha}\phi_\alpha)=\sum\ \btau_i(a_{\alpha})\phi_{\alpha}
\eqno{(3.1.1)}
$$
(Note that this extension  depends on a choice 
of a base $\{\phi_{\alpha}\}$ of the module $E$.) The fields $\{\tau_i\}$ 
form a $\Lambda E$-base of the even part $T_{\Lambda E}^{ev}$. 

We define the odd vector fields $\psi_{\alpha}\in T_{\Lambda E}^{odd}$ by 
$$
\psi_{\alpha}(\sum\ a_{\nu}\phi_{\nu})=a_{\alpha};\ \psi_{\alpha}(a)=0
\eqno{(3.1.2)}
$$
These fields form a $\Lambda E$-base of $T_{\Lambda E}^{odd}$. 

Let $\{\omega_i;\ \rho_{\alpha}\}$ be the dual base of the module 
of $1$-superforms $\Omega_{\Lambda E}=
Hom^{ev}_{\Lambda E}(T_{\Lambda E},\Lambda E)$, 
defined by 
$$
\langle\tau_i,\omega_j\rangle=\delta_{ij};\ 
\langle\psi_{\alpha},\rho_{\beta}\rangle=\delta_{\alpha\beta};\ 
\langle\tau_i,\rho_{\alpha}\rangle=
\langle\psi_{\alpha},\omega_i\rangle=0
\eqno{(3.1.3)}
$$ 

{\bf 3.2.} Let us describe the effect of a change of frame. 
Let $\fg'=\{\btau'_i;\ \phi'_{\alpha}\}$ be another frame, with 
$\btau'_i=g^{ij}\btau_j;\ \phi'_{\alpha}=A^{\alpha\beta}\phi_{\beta}$, 
$g=(g^{ij})\in GL_n(A),\ A=(A^{\alpha\beta})\in GL_m(A)$. 

The corresponding new bases $\tau'_i,$ etc. look as follows. 
$$
\tau'_i=g^{ip}\tau_p+g^{i\alpha\gamma}\phi_{\gamma}\psi_{\alpha}
\eqno{(3.2.1)}
$$
where
$$
g^{i\alpha\gamma}=g^{iq}\tau_q(A^{-1\alpha\mu})A^{\mu\gamma}
\eqno{(3.2.2)}
$$
Next,     
$$
\psi'_{\alpha}=A^{-1\mu\alpha}\psi_{\mu}
\eqno{(3.2.3)}
$$
$$
\omega'_i=g^{-1pi}\omega_p
\eqno{(3.2.4)}
$$
$$
\rho'_{\alpha}=\tau_i(A^{\alpha\gamma})\phi_{\gamma}\omega_i+
A^{\alpha\mu}\rho_{\mu}
\eqno{(3.2.5)}
$$
Formulas for the inverse transformation: 
$$
\tau_q=g^{-1qi}\tau'_i+\tau_q(A^{\alpha\gamma})\phi_{\gamma}\psi'_{\alpha}
\eqno{(3.2.6)}
$$
$$
\psi_{\beta}=A^{\alpha\beta}\psi'_{\alpha}
\eqno{(3.2.7)}
$$
$$
\omega_j=g^{pj}\omega'_p
\eqno{(3.2.8)}
$$
$$
\rho_{\beta}=A^{-1\beta\alpha}\rho'_{\alpha}+
g^{p\beta\gamma}\phi_{\gamma}\omega'_p
\eqno{(3.2.9)}
$$
These formulas show that $T$ is canonically an $A$-module quotient of 
$T_{\Lambda E}$ and  
$\Omega=Hom_A(T,A)$ is canonically an $A$-submodule 
of $\Omega_{\Lambda E}$. In fact the whole de Rham complex $\Omega^{\cdot}_A$ 
is canonically the subcomplex of  
$\Omega^{\cdot}_{\Lambda E}$. 

{\bf 3.3.} Recall that 
$$
g^{ip}\tau_p(g^{jq})=g^{jp}\tau_p(g^{iq})
\eqno{(3.3.1)}
$$
$$
g^{ip}\tau_q\tau_p(g^{jq})=g^{jq}\tau_p\tau_q(g^{ip})
\eqno{(3.3.2)}
$$
and
$$
\tau_p(g^{-1qr})=\tau_q(g^{-1pr})
\eqno{(3.3.3)}
$$
see [GII], 5.4. 

It is easy to see that 
$$
tr\bigl\{\tau_p(A)\tau_q(A^{-1})\bigr\}=
tr\bigl\{\tau_q(A)\tau_p(A^{-1})\bigr\}
\eqno{(3.3.4)}
$$
Using (3.3.1) and (3.3.4) one sees easily that 
$$ 
g^{ip}\tau_p(g^{j\nu\nu})=g^{jq}\tau_q(g^{i\nu\nu})
\eqno{(3.3.5)}
$$

{\bf 3.4.} Let $\CA=\CA_{\Lambda E;\fg}$ be the vertex superalgebroid  
corresponding to the frame $\fg$. 
 
We have the following identities in $\CA$: 

{\it $\gamma$-formulas}
 
$$
\gamma(a,b\tau_i)=-\tau_i(a)\dpar b-\tau_i(b)\dpar a
\eqno{(3.4.1)}
$$ 
$$
\gamma(a,b\phi_{\nu}\psi_{\mu})=\delta_{\nu\mu}b\dpar a
\eqno{(3.4.2)}
$$
$$
\gamma(a\phi_{\beta},b\psi_{\mu})=\delta_{\beta\mu}a\dpar b
\eqno{(3.4.3)}
$$

{\it $\langle,\rangle$-formulas} 

$$
\langle a\tau_i,b\tau_j\rangle=-b\tau_i\tau_j(a)-
a\tau_j\tau_i(b)-\tau_i(b)\tau_j(a)
\eqno{(3.4.4)}
$$
$$
\langle a\phi_{\alpha}\psi_{\beta},b\tau_i\rangle=
\delta_{\alpha,\beta}b\tau_i(a)
\eqno{(3.4.5)}
$$
$$
\langle a\phi_{\alpha}\psi_{\beta},\phi_{\alpha'}\psi_{\beta'}\rangle= 
ab\delta_{\beta\alpha'}\delta_{\beta'\alpha}
\eqno{(3.4.6)}
$$

{\it $c$-formulas} 

$$
c(a\tau_i,b\tau_j)=\frac{1}{2}\{\tau_i(b)\dpar\tau_j(a)-\tau_j(a)\dpar
\tau_i(b)\}+\frac{1}{2}\dpar\{b\tau_i\tau_j(a)-a\tau_j\tau_i(b)\}
\eqno{(3.4.7)}
$$
$$
c(a\phi_{\alpha}\psi_{\mu},b\phi_{\beta}\psi_{\nu})=
\frac{\delta_{\mu\beta}\delta_{\nu\alpha}}{2}
\bigl\{a\dpar b-b\dpar a\bigr\}
\eqno{(3.4.8)}
$$
$$
c(a\phi_{\alpha}\psi_{\mu},b\tau_i)=
-\frac{\delta_{\mu\alpha}}{2}\dpar\bigl\{b\tau_i(a)\bigr\}
\eqno{(3.4.9)}
$$
$$
c(a\tau_i,b\psi_{\alpha})=c(a\phi_{\alpha}\psi_{\mu},b\psi_{\nu})=0
\eqno{(3.4.10)}
$$  

{\bf 3.5.} Let $\fg'$ be another frame as in 3.2. We have 
$$
\gamma(a,\tau'_p)=\gamma(a,g^{pq}\tau_q+g^{p\mu\nu}\phi_{\nu}\psi_{\mu})=
$$
$$
=-\tau_q(a)\dpar g^{pq}-\tau_q(g^{pq})\dpar a+g^{p\mu\mu}\dpar a
\eqno{(3.5.1)}
$$
$$
\gamma(a\phi'_{\mu},\psi'_{\alpha})=
\gamma(aA^{\mu\beta}\phi_{\beta},A^{-1\nu\alpha}\psi_{\nu})=
aA^{\mu\beta}\dpar A^{-1\beta\alpha}
\eqno{(3.5.2)}
$$ 
Next, 
$$
\langle\tau'_i,\tau'_j\rangle=
\langle g^{ip}\tau_p+g^{i\mu\alpha}\phi_{\alpha}\psi_{\mu},
g^{jq}\tau_q+g^{j\nu\beta}\phi_{\beta}\psi_{\nu}\rangle=
$$
$$
=-2g^{ip}\tau_q\tau_p(g^{jq})-\tau_p(g^{jq})\tau_q(g^{ip})+
$$
$$
+2g^{ip}\tau_p(g^{j\nu\nu})
+g^{ip}g^{jq}\tau_p(A^{-1\mu\beta})A^{\beta\gamma}
\tau_q(A^{-1\gamma\sigma})A^{\sigma\mu}
\eqno{(3.5.3)}
$$
and
$$
\langle\tau'_i,\psi'_{\alpha}\rangle=
\langle\psi'_{\alpha},\psi'_{\beta}\rangle=0
\eqno{(3.5.4)}
$$
Finally, 
$$
c(\tau'_i,\tau'_j)=\frac{1}{2}\bigl\{\tau_p(g^{jq})\dpar\tau_q(g^{ip})-
\tau_q(g^{ip})\dpar\tau_p(g^{jq})\bigr\}+
\frac{1}{2}\bigl\{g^{i\mu\nu}\dpar g^{j\nu\mu}-
g^{j\nu\mu}\dpar g^{i\mu\nu}\bigr\}
\eqno{(3.5.5)}
$$
and
$$
c(\tau'_i,\psi'_{\alpha})=c(\psi'_{\alpha},\psi'_{\beta})=0
\eqno{(3.5.6)}
$$

\bigskip\bigskip

\centerline{\bf \S 4. Chern-Simons term}

\bigskip\bigskip

This Section is parallel to [GII], Section 5. 

{\bf 4.1.} We keep the setup of the previous section.  
Let $\CA=\CA_{\Lambda E;\fg},\ \CA'=\CA_{\Lambda E;\fg'}$ 
(resp., $\CB=\CB_{\Lambda E;\fg},\ \CB'=\CB_{\Lambda E;\fg'}$)  
be the vertex superalgebroids (resp. prealgebroids) corresponding 
to our frames. 

As in [GII], 5.5 we have a canonical isomorphism 
$$
g=g_{\fg,\fg'}=(Id_{\Lambda E},Id_{T_{\Lambda E}},Id_{\Omega_{\Lambda E}},h):\ 
\CB'\iso\CB
\eqno{(4.1.1)}
$$
where 
$$
h=h_{\fg,\fg'}:\ T_{\Lambda E}\lra\Omega_{\Lambda E}
\eqno{(4.1.2)}
$$
is defined by the condition 
$$
\langle x',h(y')\rangle=-\frac{1}{2}\langle x',y'\rangle,\ \ 
x',y'\in\{\tau'_i\}\cup\{\psi'_{\alpha}\}
\eqno{(4.1.3)}
$$
Using (3.5.3) and (3.5.4) we find the following explicit formulas for $h$: 
$$
h(\tau'_i)=h^{ij}\omega_j;\ h(\psi'_{\alpha})=0
\eqno{(4.1.4)}
$$
where $h^{ij}=h^{ij}_\Omega-h^{ij}_E$, 
$$
h^{ij}_\Omega=\tau_p\tau_j(g^{ip})+\frac{1}{2}\tau_q(g^{ip})\tau_p(g^{rq})g^{-1j
r}
\eqno{(4.1.5)}
$$
cf. [GII], (5.7.2), and 
$$
h^{ij}_E= \tau_j(g^{i\nu\nu}) +
\frac{1}{2}g^{iq}\tau_j(A^{-1\mu\beta})A^{\beta\gamma}
\tau_q(A^{-1\gamma\nu})A^{\nu\mu}
\eqno{(4.1.6)}
$$
The meaning of the notation $h_{\Omega}$ will become clear below, 
see \S 6.   

{\bf 4.2.} We have 
$$
\CA=g_*\CA'\dplus\fb
\eqno{(4.2.1)}
$$
where the closed $3$-form $\fb\in\Omega^{3,cl}_{\Lambda E}$ is defined 
by
$$
\fb(x',y')=c(x',y')-x'(h(y'))+(-1)^{p(x')p(y')}y'(h(x')),
\eqno{(4.2.2)}
$$ 
$x',y'\in\{\omega'_i\}\cup\{\psi'_{\alpha}\}$, cf. [GII], (5.7.3). 
 
It is easy to see that $\psi'_{\alpha}(h(\tau'_i))=0$; 
on the other hand we know already that $h(\psi'_{\alpha})=0$ and 
$c(\psi'_{\alpha},y')=0$. It follows that $\fb\in\Omega^{3,cl}\subset 
\Omega^{3,cl}_{\Lambda E}$. 

Next, we have 
$$
\tau'_i(h(\tau'_j))=(g^{ip}\tau_p+
g^{i\alpha\gamma}\phi_{\gamma}\psi_{\alpha})(h^{jq}\omega_q)=
$$
(note that the second summand is zero)
$$
=(g^{ip}\tau_p)(h_{\Omega}^{jq}\omega_q)+(g^{ip}\tau_p)(h_E^{jq}\omega_q)=
g^{ip}\tau_p(h^{jq}_{\Omega})\omega_q+h^{jp}_{\Omega}\dpar g^{ip}-
$$
$$
-g^{ip}\tau_p(h^{jq}_E)\omega_q-h^{jp}_E\dpar g^{ip}
\eqno{(4.2.3)}
$$
It follows that $\fb=\fb_{\Omega}-\fb_E$ where $\fb_{\Omega}, \fb_E\in\Omega^3$ 
are 
given by 
$$
\fb_{\Omega}(\tau'_i,\tau'_j)=\frac{1}{2}\bigl\{\tau_p(g^{jq})\dpar\tau_q(g^{ip}
)-
\tau_q(g^{ip})\dpar\tau_p(g^{jq})\bigr\}-
$$
$$ 
-g^{ip}\tau_p(h^{jq}_{\Omega})\omega_q-h^{jp}_{\Omega}\dpar g^{ip}+
g^{jp}\tau_p(h^{iq}_{\Omega})\omega_q+h^{ip}_{\Omega}\dpar g^{jp}
\eqno{(4.2.4)}
$$
and
$$
\fb_E(\tau'_i,\tau'_j)=
-\frac{1}{2}\bigl\{g^{i\mu\nu}\dpar g^{j\nu\mu}-
g^{j\nu\mu}\dpar g^{i\mu\nu}\bigr\}-
$$
$$
-g^{ip}\tau_p(h^{jq}_E)\omega_q
+g^{jp}\tau_p(h^{iq}_E)\omega_q
-h^{jp}_E\dpar g^{ip}+h^{ip}_E\dpar g^{jp}
\eqno{(4.2.5)}
$$   
The form $\fb_{\Omega}$ has already been computed in [GII], Magic Lemma 5.6 and 
Theorem 6.4 (b),  
and is equal to 
$$
\fb_{\Omega}(\tau'_i,\tau'_j)=-\frac{1}{2}tr\bigl\{
g^{-1}\tau'_i(g)g^{-1}\tau'_j(g)g^{-1}\tau'_r(g)-
g^{-1}\tau'_j(g)g^{-1}\tau'_i(g)g^{-1}\tau'_r(g)\bigr\}\omega'_r
\eqno{(4.2.6)}
$$
cf. {\it loc.cit.} (5.5.3) and (6.4.2). Note that $\fb_{\Omega}$ is closed, 
hence 
$\fb_E$ is closed.

{\bf 4.3. Magic Lemma.} {\it We have
$$
\fb_E(\tau'_i,\tau'_j)=-\frac{1}{2}tr\bigl\{
A^{-1}\tau'_i(A)A^{-1}\tau'_j(A)A^{-1}\tau'_r(A)-
A^{-1}\tau'_j(A)A^{-1}\tau'_i(A)A^{-1}\tau'_r(A)\bigr\}\omega'_r
\eqno{(4.3.1)}
$$}

{\it Proof.} Let us denote the six terms in (4.2.5) by $A, A', B, B', C$ and 
$C'$. We have 
$$
A=-\frac{1}{2}g^{iq}\tau_q(A^{-1\mu\alpha})A^{\alpha\nu}\tau_r\bigl\{
g^{jp}\tau_p(A^{-1\nu\beta})A^{\beta\mu}\bigr\}\omega_r=
$$
$$
=-\frac{1}{2}g^{iq}\tau_q(A^{-1\mu\alpha})A^{\alpha\nu}\bigl\{
\tau_r(g^{jp})\tau_p(A^{-1\nu\beta})A^{\beta\mu}\omega_r+
g^{jp}\tau_r\tau_p(A^{-1\nu\beta})A^{\beta\mu}\omega_r+
$$
$$
+g^{jp}\tau_p(A^{-1\nu\beta})\tau_r(A^{\beta\mu})\omega_r\bigr\}
$$
Next, 
$$
B= -\frac{1}{2}g^{ip}\tau_p\bigl\{
g^{jq}\tau_r(A^{-1\mu\beta})A^{\beta\gamma}\tau_q(A^{-1\gamma\nu})A^{\nu\mu}
\bigr\}\omega_r - g^{ip}\tau_p\tau_r(g^{j\nu\nu})\omega_r=
$$
$$
=-\frac{1}{2}g^{ip}\tau_p(g^{jq})
\tau_r(A^{-1\mu\beta})A^{\beta\gamma}\tau_q(A^{-1\gamma\nu})A^{\nu\mu}
\omega_r
-\frac{1}{2}g^{ip}g^{jq}
\tau_p\tau_r(A^{-1\mu\beta})A^{\beta\gamma}\tau_q(A^{-1\gamma\nu})A^{\nu\mu}
\omega_r-
$$
$$   
-\frac{1}{2}g^{ip}g^{jq}
\tau_r(A^{-1\mu\beta})\tau_p(A^{\beta\gamma})\tau_q(A^{-1\gamma\nu})A^{\nu\mu}
\omega_r-
-\frac{1}{2}g^{ip}g^{jq}
\tau_r(A^{-1\mu\beta})A^{\beta\gamma}\tau_p\tau_q(A^{-1\gamma\nu})A^{\nu\mu}
\omega_r-
$$
$$
-\frac{1}{2}g^{ip}g^{jq}
\tau_r(A^{-1\mu\beta})A^{\beta\gamma}\tau_q(A^{-1\gamma\nu})
\tau_p(A^{\nu\mu})\omega_r - g^{ip}\tau_p\tau_r(g^{j\nu\nu})\omega_r
$$
Finally, 
$$
C= -\frac{1}{2}g^{jq}\tau_p(A^{-1\mu\beta})A^{\beta\gamma}
\tau_q(A^{-1\gamma\nu})A^{\nu\mu}\tau_r(g^{ip})\omega_r - 
\tau_r(g^{ip})\tau_p(g^{j\nu\nu})\omega_r
$$
We see that $A1=-C1', A2=-B'2, B1=-B'1$ (by (3.4.6)), $B2=-A'2, B4=-B'4$ and 
$C1=-A'1$. Next, 
$$
B6+C2=-\tau_r\bigl\{g^{ip}\tau_p(g^{j\nu\nu})\bigr\}, 
$$
so $B6+C2=-B'6-C'2$ by (3.3.5).  

So we are left with three terms: $A3, B3$ and $B5$ and 
their primed partners. It is easy to see that
$$
A3=B3=-B5=-\frac{1}{2}tr\bigl\{
A^{-1}\tau'_i(A)A^{-1}\tau'_j(A)A^{-1}\tau'_r(A)\bigr\}\omega'_r,
$$
which implies the Lemma. $\btu$

\bigskip\bigskip

\centerline{\bf \S 5. Atiyah term}

\bigskip\bigskip

This section is parallel to [GII], Section 6. 

{\bf 5.1.} We keep the setup of the previous section. 
Let us denote by $\CT_{\Lambda E}=(\Lambda E, T_{\Lambda E}, 
\Omega_{\Lambda E}, \dpar)$ the extended vertex superalgebroid Lie 
corresponding to our data. 

Let $\fg''=\{\btau''_i,\phi''_{\alpha}\}$ be a third frame of $(A,E)$, 
with $\btau''_i=g^{\prime ij}\btau'_j, \phi''_{\alpha}=A^{\prime\alpha\beta}
\phi'_{\beta}$. We have the corresponding new bases of $T_{\Lambda E}$ and 
$\Omega_{\Lambda E}$ given by
$$
\tau''_i=g^{\prime ip}\tau'_p+g^{\prime i\alpha\gamma}
\phi'_{\gamma}\psi'_{\alpha}
\eqno{(5.1.1)}
$$
where
$$
g^{\prime i\alpha\gamma}=g^{\prime iq}\tau'_q(A^{\prime -1\alpha\mu})
A^{\prime\mu\gamma}
\eqno{(5.1.2)}
$$     
$$
\psi''_{\alpha}=A^{\prime -1\mu\alpha}\psi'_{\mu}
\eqno{(5.1.3)}
$$
$$
\omega''_i=g^{\prime -1pi}\omega'_p
\eqno{(5.1.4)}
$$
$$
\rho''_{\alpha}=\tau'_i(A^{\prime\alpha\gamma})\phi'_{\gamma}\omega'_i+
A^{\prime\alpha\mu}\rho'_{\mu}
\eqno{(5.1.5)}
$$

{\bf 5.2.} Let $\CA''=\CA_{\Lambda E,\fg''}$ (resp. 
$\CB''=\CB_{\Lambda E;\fg''}$) be the vertex superalgebroid (resp. 
prealgebroid) corresponding to the third frame. We have 
canonical isomorphisms 
$$
\CB''\buildrel{g_{\fg',\fg''}}\over\iso \CB'
\buildrel{g_{\fg,\fg'}}\over\iso \CB
$$
as well as the morphism $g_{\fg,\fg''}:\ \CB''\iso\CB$ over 
$Id_{\CT_{\Lambda E}}$, given by functions 
$h_{\fg,\fg'}, h_{\fg',\fg''}, h_{\fg,\fg''}$, 
and we are aiming to compute the discrepancy 
$$
\fa=\fa_{\fg,\fg',\fg''}:= h_{\fg,\fg'}+h_{\fg',\fg''}-h_{\fg,\fg''}
\in\Omega^2_{\Lambda E}
\eqno{(5.2.1)}
$$
Note that 
$$
\gamma(a,b\psi_{\mu})=-\psi_{\mu}(a)\dpar b-\psi_{\mu}(b)\dpar a=0
\eqno{(5.2.2)}
$$
(in $\CA$), hence
$$
h_{\fg,\fg'}(a\psi'_{\mu})=ah_{\fg,\fg'}(\psi'_{\mu})-
\gamma(a,\psi'_{\mu})=0
\eqno{(5.2.3)}
$$
therefore
$$
h_{\fg,\fg'}(\psi''_{\alpha})=h_{\fg,\fg'}(A^{\prime -1\mu\alpha}\psi'_{\mu})
=0
\eqno{(5.2.4)}
$$
It follows that 
$$
\fa(\psi''_{\alpha})=0
\eqno{(5.2.5)}
$$ 

{\bf 5.3.}  Let us denote for 
brevity $h:=h_{\fg,\fg'},\ h'=h_{\fg',\fg''},\ h''=h_{\fg,\fg''}$. 

We have
$$
h(a\tau'_p)=ah(\tau'_p)-\gamma(a,\tau'_p)=
$$
$$
=ah^{pr}\omega_r+
\tau_q(a)\dpar g^{pq}+\tau_q(g^{pq})\dpar a-g^{p\mu\mu}\dpar a
\eqno{(5.3.1)}
$$
and
$$
h(a\phi'_{\gamma}\psi'_{\alpha})=-aA^{\gamma\beta}\dpar A^{-1\beta\alpha}
\eqno{(5.3.2)}
$$
Thus we get
$$
h(\tau''_i)=h(g^{\prime ip}\tau'_p+g^{\prime i\alpha\gamma}\phi'_{\gamma}
\psi'_{\alpha})=
$$
$$
=g^{\prime ip}h^{pr}_{\Omega}\omega_r
-g^{\prime ip}h^{pr}_E\omega_r+\tau_q(g^{\prime ip})\dpar g^{pq}+
\tau_q(g^{pq})\dpar g^{\prime ip}-
$$
$$
-g^{p\mu\mu}\dpar g^{\prime ip}-
g^{\prime i\alpha\gamma}A^{\gamma\beta}\dpar A^{-1\beta\alpha}
\eqno{(5.3.3)}
$$
It follows that 
$$
\fa=\fa_{\Omega}-\fa_E
\eqno{(5.3.4)}
$$
where
$$
\fa_{\Omega}(\tau''_i)=
g^{\prime ip}h^{pr}_{\Omega}\omega_r 
+\tau_q(g^{\prime ip})\dpar g^{pq}+
\tau_q(g^{pq})\dpar g^{\prime ip}+
h_{\Omega}^{\prime ip}\omega'_p-h_{\Omega}^{\prime\prime ir}\omega_r
\eqno{(5.3.5)}
$$
and
$$
\fa_E(\tau''_i)=g^{\prime ip}h^{pr}_E\omega_r
+g^{p\mu\mu}\dpar g^{\prime ip}+
g^{\prime i\alpha\gamma}A^{\gamma\beta}\dpar A^{-1\beta\alpha}
+h_E^{\prime ip}\omega'_p-h_E^{\prime\prime ir}\omega_r
\eqno{(5.3.6)}
$$
The form $\fa_{\Omega}$ has already been computed in [GII], Section 6. Namely, 
by Theorem 6.4 (a) from {\it op. cit.}, 
$$
\fa_{\Omega}(\tau''_i)=\frac{1}{2}tr\bigl\{
g^{\prime -1}\tau''_i(g')\tau''_s(g)g^{-1}- 
g^{\prime -1}\tau''_s(g')\tau''_i(g)g^{-1}\bigl\}\omega''_s
\eqno{(5.3.7)}
$$

{\bf 5.4.} Let us compute $\fa_E(\tau''_i)$. Let us denote the five terms 
in the rhs of (5.3.6) by $\fA, \fB, \fC, \fD$ and $\fE$. Thus, 
$$
\fA=g^{\prime ip}h^{pr}_E\omega_r=\frac{1}{2}g^{\prime ip}g^{pq}
\tau_r(A^{-1\mu\beta})A^{\beta\gamma}\tau_q(A^{-1\gamma\nu})
A^{\nu\mu}\omega_r+g^{\prime ip}\tau_r(g^{p\nu\nu})\omega_r
\eqno{(5.4.1)}
$$
$$
\fB=g^{p\mu\mu}\dpar g^{\prime ip}=
g^{pq}\tau_q(A^{-1\mu\alpha})A^{\alpha\mu}\tau_r(g^{\prime ip})\omega_r
\eqno{(5.4.2)}
$$
$$
\fC=g^{\prime i\alpha\gamma}A^{\gamma\beta}\dpar A^{-1\beta\alpha}=
g^{\prime ip}\tau'_p(A^{\prime -1\alpha\mu})A^{\prime\mu\gamma}
A^{\gamma\beta}\tau_r(A^{-1\beta\alpha})\omega_r
\eqno{(5.4.3)}
$$
$$
\fD=h_E^{\prime ij}\omega'_j=
\frac{1}{2}g^{\prime is}\tau'_j(A^{\prime -1\mu\beta})A^{\prime\beta\gamma}
\tau'_s(A^{\prime -1\gamma\nu})A^{\prime\nu\mu}g^{-1rj}\omega_r+
\tau'_j(g^{\prime i\nu\nu})g^{-1rj}\omega_r=
$$
$$
=\frac{1}{2}g^{\prime is}\tau_r(A^{\prime -1\mu\beta})
A^{\prime\beta\gamma}g^{sq}\tau_q(A^{\prime -1\gamma\nu})
A^{\prime\nu\mu}\omega_r+
\tau_r\bigl\{g^{\prime ip}g^{pq}\tau_q(A^{\prime -1\nu\mu})A^{\mu\nu}\bigr\}
\omega_r
\eqno{(5.4.4)}
$$
and
$$
\fE=-h^{\prime\prime ir}\omega_r=-\frac{1}{2}(g'g)^{iq}
\tau_r\bigl((A'A)^{-1\mu\beta}\bigr)(A'A)^{\beta\gamma}
\tau_q\bigl((A'A)^{-1\gamma\nu}\bigr)(A'A)^{\nu\mu}\omega_r 
-\tau_r(g^{\prime\prime i\nu\nu})=
$$
$$
=-\frac{1}{2}(g'g)^{iq}\bigl\{
\tau_r(A^{-1\mu\sigma})A^{\prime -1\sigma\beta}+
A^{-1\mu\sigma}\tau_r(A^{\prime -1\sigma\beta})\bigr\}
A^{\prime\beta\rho}A^{\rho\gamma}\times
$$
$$
\times\bigl\{
\tau_q(A^{-1\gamma\delta})A^{\prime -1\delta\nu}+
A^{-1\gamma\delta}\tau_q(A^{\prime -1\delta\nu})\bigr\}
A^{\prime\nu\epsilon}A^{\epsilon\mu}\omega_r 
-\tau_r(g^{\prime\prime i\nu\nu})=
$$
$$
=-\frac{1}{2}g^{\prime ip}g^{pq}\bigl\{
\tau_r(A^{-1\mu\sigma})A^{\sigma\gamma}\tau_q(A^{-1\gamma\delta})
A^{\delta\mu}\omega_r+
\tau_r(A^{-1\mu\sigma})\tau_q(A^{\prime -1\sigma\nu})
A^{\prime\nu\epsilon}A^{\epsilon\mu}\omega_r+
$$
$$
+\tau_r(A^{\prime -1\sigma\beta})A^{\prime\beta\rho}
A^{\rho\gamma}\tau_q(A^{-1\gamma\sigma})\omega_r+
\tau_r(A^{\prime -1\sigma\beta})A^{\prime\beta\rho}
\tau_q(A^{\prime -1\rho\nu})A^{\prime\nu\sigma}\omega_r\bigr\}-
$$
$$
-\tau_r\bigl\{(g'g)^{iq}\tau_q\bigl((A'A)^{-1\nu\mu}\bigr)
(A'A)^{\mu\nu}\bigr\}\omega_r
\eqno{(5.4.5)}
$$
We see that $\fA1=-\fE1,\ \fC=-2\fE2,\ \fD1=-\fE4$. 
It is easy to see that $\fA2+\fD2+\fE5=-\fB$. 
 
Finally, 
$$
\fC+\fE2=\frac{1}{2}tr\bigl\{
A^{\prime -1}\tau''_i(A')\tau''_r(A)A^{-1}\bigr\}\omega''_r
\eqno{(5.4.6)}
$$
and
$$
\fE3=-\frac{1}{2}tr\bigl\{
A^{\prime -1}\tau''_r(A')\tau''_i(A)A^{-1}\bigr\}\omega''_r
\eqno{(5.4.7)}
$$
So, we have proven 

{\bf 5.5. Lemma.} {\it The form $\fa_E$ is given by  
$$
\fa_E(\tau''_i)=\frac{1}{2}tr\bigl\{
A^{\prime -1}\tau''_i(A')\tau''_r(A)A^{-1} 
- A^{\prime -1}\tau''_r(A')\tau''_i(A)A^{-1}\bigr\}\omega''_r
\eqno{(5.5.1)}
$$}

Combining 4.3 and 5.5 we get

{\bf 5.6. Theorem.} (a) {\it The cocycle $\fa_{\fg,\fg',\fg''}$ is given by 
$$
\fa_{\fg,\fg',\fg''}(\tau''_i)=\frac{1}{2}tr\bigl\{
g^{\prime -1}\tau''_i(g')\tau''_r(g)g^{-1}- 
g^{\prime -1}\tau''_r(g')\tau''_i(g)g^{-1}\bigl\}\omega''_r - 
$$
$$
-\frac{1}{2}tr\bigl\{
A^{\prime -1}\tau''_i(A')\tau''_r(A)A^{-1} 
- A^{\prime -1}\tau''_r(A')\tau''_i(A)A^{-1}\bigr\}\omega''_r
\eqno{(5.6.1)}
$$} 

(b) {\it The $3$-form $\fb_{\fg,\fg'}$ is given by 
$$
\fb_{\fg,\fg'}(\tau'_i,\tau'_j)=-\frac{1}{2}tr\bigl\{
g^{-1}\tau'_i(g)g^{-1}\tau'_j(g)g^{-1}\tau'_r(g)-
g^{-1}\tau'_j(g)g^{-1}\tau'_i(g)g^{-1}\tau'_r(g)\bigr\}\omega'_r + 
$$
$$
+\frac{1}{2}tr\bigl\{
A^{-1}\tau'_i(A)A^{-1}\tau'_j(A)A^{-1}\tau'_r(A)-
A^{-1}\tau'_j(A)A^{-1}\tau'_i(A)A^{-1}\tau'_r(A)\bigr\}\omega'_r
\eqno{(5.6.2)}
$$}

{\bf 5.7. Lemma.} {\it Let $E^*=Hom_A(E,A)$ be the dual module. We have 
$$
(\fa_{E^*},\fb_{E^*})=(\fa_E,\fb_E)
\eqno{(5.7.1)}
$$}

Indeed, this follows from the easy identities
$$
tr\bigl\{A^t\tau_i((A^t)^{-1})A^t\tau_j((A^t)^{-1})A^t\tau_r((A^t)^{-1})
\bigr\}=
-tr\bigl\{A^{-1}\tau_r(A)A^{-1}\tau_j(A)A^{-1}\tau_i(A)\bigr\}
\eqno{(5.7.2)}
$$
and
$$
tr\bigl\{A^t\tau_i((A^t)^{-1}\tau_j((B^t)^{-1})B^t\bigr\}=
tr\bigl\{A^{-1}\tau_i(A)\tau_j(B)B^{-1}\bigr\}
\eqno{(5.7.3)}
$$

{\bf 5.8.} Let us pass to the global situation. Let $X$ be a 
smooth variety over $k$ and $E$ be a vector bundle over $X$. 
As in [GII], we define the gerbe $\fD_{\Lambda E}$ of chiral 
differential operators on $\Lambda E$ over $X$. 

Its characteristic class $c(\fD_{\Lambda E})$ will belong 
to the second hypercohomology $H^2(X;\Omega^{[2,3\rangle}_{\Lambda E})$ 
(in obvious notations). Recall that we have a canonical 
imbedding of de Rham complexes
$$
\Omega^{\cdot}_X\hra \Omega^{\cdot}_{\Lambda E}
\eqno{(5.8.1)}
$$
In [GII], 7.6 we have defined the "Atiyah-Chern-Simons" characteristic class 
$c(E)\in H^2(X;\Omega^{[2,3\rangle}_X)$; let us denote by 
$c(E)_{\Lambda E}$ its image in 
$H^2(X;\Omega^{[2,3\rangle}_{\Lambda E})$. 

The theorem below is an immediate consequence Theorem 5.6 and Lemma 5.7. 

{\bf 5.9. Theorem.} {\it The class $c(\fD_{\Lambda E})$ is equal to 
$$
c(\fD_{\Lambda E})=c(\Theta_X)-c(E)=c(\Omega^1_X)-c(E)
\eqno{(5.9.1)}
$$
where $\Theta_X$ is the tangent bundle.}

\bigskip\bigskip

\newpage

\centerline{\bf \S 6. Chiral de Rham complex}

\bigskip\bigskip

{\bf 6.1.} Let us return to the local situation 3.1, 4.1.  Let  
$E$ be equal to the module of vector fields $T$. 
Given a base $\{\btau_i\}$ consisting of commuting vector fields, 
we get a frame $\fg=\{\btau_i; \phi_i:=\btau_i\}$ of $(A,E)$. Let us call such 
frames {\it natural}. 

Let $\fg, \fg'$ be two natural frames, with transition matrices as in 
3.2. By definition, $(A^{rs})=(g^{rs})$. Therefore the coefficients 
$g^{i\alpha\gamma}$ (3.2.2) are given by  
$$
g^{i\alpha\gamma}=g^{iq}\tau_q(g^{-1\alpha\mu})g^{\mu\gamma}=
g^{iq}\tau_{\alpha}(g^{-1q\mu})g^{\mu\gamma}=
$$
$$
=-\tau_{\alpha}(g^{iq})g^{-1q\mu}g^{\mu\gamma}=
-\tau_{\alpha}(g^{i\gamma})
\eqno{(6.1.1)}
$$
where we have used (3.3.3). Consequently the function $h_E$ (4.1.6) 
is given by  
$$
h_E^{ij}=-\tau_j\tau_{\nu}(g^{i\nu})+
\frac{1}{2}g^{iq}\tau_j(g^{-1\mu\beta})g^{\beta\gamma}
\tau_q(g^{-1\gamma\nu})g^{\nu\mu} 
\eqno{(6.1.2)}
$$
The second summand is equal to 
$$
-\frac{1}{2}g^{iq}\tau_j(g^{-1\mu\beta})\tau_q(g^{\beta\gamma})
g^{-1\gamma\nu}g^{\nu\mu}=
-\frac{1}{2}g^{iq}\tau_j(g^{-1\mu\beta})\tau_q(g^{\beta\mu})=
$$
$$
=-\frac{1}{2}g^{iq}\tau_\mu(g^{-1j\beta})\tau_q(g^{\beta\mu})=
-\frac{1}{2}\tau_{\mu}(g^{-1j\beta})g^{\beta q}\tau_q(g^{i\mu})=
\frac{1}{2}g^{-1j\beta}\tau_{\mu}(g^{\beta q})\tau_q(g^{i\mu}) 
\eqno{(6.1.3)}
$$ 
We see that the first summand in (6.1.2) is equal to minus the first summand 
of (4.1.5), and second summand of (6.1.2) is equal to the second summand 
of (4.1.5). Thus  
$$
h^{ij}=2\tau_p\tau_j(g^{ip})
\eqno{(6.1.4)}
$$  
If $\fg, \fg', \fg''$ are natural frames of $(A,T)$ then 
Theorem 5.6 says that $\fa_{\fg,\fg',\fg''}=\fb_{\fg,\fg'}=0$. 
This means that 

{\it the chiral superalgebroid $\CA_{\Lambda T;\fg}$ 
does not depend, up to a canonical isomorphism, on the choice of the base 
$\{\btau_i\}$. In other words, we have a canonically defined chiral algebroid 
$\CA_{\Lambda T}$.} 

Passing to chiral envelopes, we get a canonically defined 
chiral (vertex) superalgebra $D^{ch}_{\Lambda T}$ of {\it chiral 
differential operators on $T$}.  

It follows that 

{\it for each smooth variety $X$ we have a canonically defined 
sheaf of chiral superalgebras $\CD^{ch}_{\Lambda \Theta_X}$.  
These Zariski sheaves form in fact a sheaf in the \'etale topology.} 

The gluing functions for this sheaf are given explicitly by (6.1.4). 

{\bf 6.2. Lemma.} {\it In the situation} 4.1, {\it consider the function 
$h_E$ 
$$
h^{ij}_E= \tau_j(g^{i\nu\nu}) +
\frac{1}{2}g^{iq}\tau_j(A^{-1\mu\beta})A^{\beta\gamma}
\tau_q(A^{-1\gamma\nu})A^{\nu\mu}
\eqno{(6.2.1)}
$$
The function $h_{E^*}$ associated with the dual module $E^*$ is given by 
$$
h^{ij}_{E^*}= -\tau_j(g^{i\nu\nu}) +
\frac{1}{2}g^{iq}\tau_j(A^{-1\mu\beta})A^{\beta\gamma}
\tau_q(A^{-1\gamma\nu})A^{\nu\mu}
\eqno{(6.2.2)}
$$}

This follows from the identities 
$$
tr\bigl\{\tau_i(A^t)A^{-1}\bigr\}=
-tr\bigl\{\tau_i(A^{-1})A\bigr\}
\eqno{(6.2.3)}
$$
and
$$
tr\bigl\{\tau_i(A^t)A^{-1}\tau_j(A^t)A^{-1}\bigr\}= 
tr\bigl\{\tau_i(A^{-1})A\tau_j(A^{-1})A\bigr\}
\eqno{(6.2.4)}
$$
(cf. (5.7.2), (5.7.3)).  
     
{\bf 6.3.} Let $E$ be the module of $1$-forms $\Omega=\Omega^1_{A/k}$;  
its exterior algebra 
is the de Rham algebra of differential forms 
$\Omega^{\cdot}=\Omega^{\cdot}_{A/k}$. Frames of the form 
$\fg=\{\btau_i,\phi_i:=\omega_i\}$ will be called natural.  

If $\fg, \fg'$ are natural frames then 
formulas (6.1.2) and (6.1.3), together with the previous lemma, 
show that $h^{ij}_E=h^{ij}_{\Omega}$ where $h_{\Omega}$ is given by (4.1.5). 
(Of course one easily checks this directly.) This explains the notation 
for $h^{ij}_{\Omega}$). 

In other words, we arrive at an interesting conclusion. 

{\bf 6.4. Theorem.} {\it The matrices $h=(h^{ij})$ defined in} 4.1 {\it 
are equal to $0$ if $E=\Omega$ and frames $\fg, \fg'$ are natural.} 

{\bf 6.4.1.} {\it Warning.} The functions $h_{\fg,\fg'}$ are 
{\it non}zero since they are not linear. 

{\bf 6.5.} On the other hand, Theorem 5.6 together Lemma 5.7 say that 
$\fa_E$ and $\fb_E$ are $0$ for $E=\Omega$ (for natural frames $\fg, \fg', 
\fg''$ of $(A,\Omega)$). This gives us a {\it canonically defined} 
chiral superalgebroid $\CA_{\Omega^{\cdot}}$. Its vertex envelope 
will be denoted $D^{ch}_{\Omega^{\cdot}}$ and called {\it the chiral 
algebra of differential operators on $\Omega^{\cdot}$.}   

This implies 

{\bf 6.6. Theorem.} {\it For each smooth variety $X$ the construction} 
6.2 - 6.5 {\it gives a canonically defined 
sheaf of chiral superalgebras $\CD^{ch}_{\Omega^{\cdot}_X}$. 
  These Zariski sheaves form a sheaf in the \'etale topology.} 

{\bf 6.7.} The de Rham differential may be considered as an 
odd first order differential operator acting on $\Omega^{\cdot}_X$ 
(and commuting with itself). 

In coordinates, if $\fg=\{\btau_i;\phi_i:=\omega_i\}$ is a natural 
frame of $(A,\Omega)$, it is given by 
$$
Q^{cl}_{\fg}=\phi_i\tau_i
\eqno{(6.7.1)}
$$
($cl$ is for "classical"). 
Let us check the independence of (6.7.1) on the choice of a frame. 
Let $\fg'=\{\btau'_i;\omega'_i\}$ be another natural frame, 
$\btau'_i=g^{ij}\btau_i;\ \omega'_i=g^{-1pi}\omega_p$ (cf. (3.2.4)). 

Using (3.2.1) we have
$$
Q^{cl}_{\fg'}=\phi'_i\tau'_i=g^{-1pi}\phi_p\bigl\{
g^{iq}\tau_q+g^{irs}\phi_s\psi_r\bigl\}
\eqno{(6.7.2)}
$$
where
$$
g^{ipq}=g^{ir}\tau_r(g^{sp})g^{-1qs}=
g^{sr}\tau_r(g^{ip})g^{-1qs}=
\tau_q(g^{ip})
\eqno{(6.7.3)}
$$
(we have used (3.3.1)). So, 
$$
Q^{cl}_{\fg'}=\phi_p\tau_p+g^{-1pi}g^{irs}\phi_s\psi_r=
Q^{cl}_{\fg}+g^{-1pi}\tau_s(g^{ir})\phi_p\phi_s\psi_r
\eqno{(6.7.4)}
$$
Note that 
$$
g^{-1pi}\tau_s(g^{ir})=-\tau_s(g^{-1pi})g^{ir}
$$
which is symmetric under the permutation of $p$ with $s$, due to (3.3.3); 
therefore the second summand in (6.7.4) is zero, 
i.e. $Q^{cl}_{\fg}=Q{cl}_{\fg'}$. 

Thus, $Q^{cl}$ is a correctly defined odd element of $D_{\Omega^{\cdot}}$. 
It is obvious from (6.7.1) that $[Q^{cl},Q^{cl}]=0$. 

{\bf 6.8.} Let us investigate the chiral counterpart of $Q^{cl}$. 
Let us define an odd element $Q_{\fg}$ (of conformal weight $1$) 
of the vertex superalgebra 
$D^{ch}_{\Omega^{\cdot};\fg}:=U\CA_{\Omega^{\cdot};\fg}$ by 
$$
Q_{\fg}=\phi_{i(-1)}\tau_i
\eqno{(6.8.1)}
$$
Let $\fg'$ be another natural frame as in 6.7. Due to 
Theorem 6.4, the element $Q_{\fg'}$ goes  under the canonical 
isomorphism $D^{ch}_{\Omega^{\cdot};\fg'}=D^{ch}_{\Omega^{\cdot};\fg}$   
to 
$$
Q_{\fg'}=\phi'_{i(-1)}\tau'_i=\phi'_i\tau'_i-\gamma(\phi'_i,\tau'_i)
\eqno{(6.8.2)}
$$
cf. [GII], (3.3.1). 

{\bf 6.9. Lemma.} {\it We have (in $D_{\Omega^{\cdot};\fg}$)  
$$
\gamma(\phi'_i,\tau'_i)=-\dpar\bigl\{
tr\bigl(\tau_r(g)g^{-1}\bigr)\phi_r\bigr\}
\eqno{(6.9.1)}
$$}

{\bf 6.10.} Before the proof, let us write down useful formulas
$$
\gamma(a\phi_r,b\tau_i)=-\tau_i(a)\phi_r\dpar b-\tau_i(b)\dpar(a\phi_r)
\eqno{(6.10.1)}
$$
and 
$$
\gamma(a\phi_r,b\phi_s\psi_p)=-\delta_{rp}a\dpar(b\phi_s)+
\delta_{sp}b\dpar(a\phi_r)
\eqno{(6.10.2)}
$$

{\bf 6.11.} {\it Proof of} 6.9. We have 
$$
\gamma(\phi'_i,\tau'_i)=\gamma(g^{-1qi}\phi_q,
g^{ip}\tau_p+g^{isr}\phi_r\psi_s)
$$   
where
$$
\gamma(g^{-1qi}\phi_q,g^{ip}\tau_p)=
-\tau_p(g^{-1qi})\phi_q\dpar g^{ip}-\tau_p(g^{ip})\dpar(g^{-1qi}\phi_q)
\eqno{(6.11.1)}
$$
and
$$
\gamma(g^{-1qi}\phi_q,g^{isr}\phi_r\psi_s)=
-g^{-1qi}\dpar(g^{iqr}\phi_r)+
g^{irr}\dpar(g^{-1qi})\phi_q
\eqno{(6.11.2)}
$$
Since $g^{irr}=\tau_r(g^{ir})$, the second summands in (6.11.1) and 
(6.11.2) cancel out. On the other hand, the first term in (6.11.1) 
is equal to 
$$
-\tau_q(g^{-1ri})\tau_s(g^{iq})\omega_s\phi_r= 
-\tau_r(g^{-1si})\tau_s(g^{iq})\omega_s\phi_r=
$$
$$
=g^{-1qa}\tau_r(g^{ab})g^{-1bi}\tau_s(g^{iq})\omega_s\phi_r= 
-\tau_r(g^{ab})\tau_s(g^{-1ba})\omega_s\phi_r=
-\tau_r(g^{ab})\dpar(g^{-1ba})\phi_r
$$
Therefore
$$
\gamma(\phi'_i,\tau'_i)=-\tau_r(g^{ab})\dpar(g^{-1ba})\phi_r-
g^{-1ba}\dpar\bigl\{\tau_r(g^{ab})\phi_r\bigr\}=
-\dpar\bigl\{\tau_r(g^{ab})g^{-1ba}\phi_r\bigr\}, 
$$
QED. 

{\bf 6.12.} From (6.8.1) we have $Q_{\fg}=\phi_i\tau_i$, and 
from 6.7 $\phi_i'\tau_i'=\phi_i\tau_i$. Therefore, (6.8.2) and 
Lemma 6.9 imply 

{\bf 6.13. Theorem.} {\it We have 
$$
Q_{\fg'}=Q_{\fg}+\dpar\bigl\{
tr\bigl(\tau_r(g)g^{-1}\bigr)\phi_r\bigr\}
\eqno{(6.13.1)}
$$}

{\bf 6.14.} Consider the field $Q_{\fg}(z)$ acting on the vertex algebra 
$D^{ch}_{\Omega^{\cdot}}$. Due to (6.13.1), its zeroth component 
$Q_{\fg 0}$ does not depend on the choice of the frame $\fg$. 
Therefore we get a canonical operator $Q_0$ acting on 
$D^{ch}_{\Omega^{\cdot}}$. 

Since it is a zeroth component of a field, it is a derivation 
of the vertex algebra, and it is obvious from the local definition 
(6.8.1) that $[Q_0,Q_0]=0$. 

Consequently, for each smooth variety $X$ we get a canonical odd derivation 
$Q_{0X}$ of the sheaf $\CD^{ch}_{\Omega^{\cdot}_X}$, such that 
$[Q_{0X},Q_{0X}]=0$.  
The pair $(\CD^{ch}_{\Omega^{\cdot}_X},Q_{0X})$ is the 
{\it chiral de Rham complex} 
from [MSV]. 

Our Theorem 6.13 is a version of {\it op. cit.}, (4.1c). 

{\bf 6.15.} In the situation 6.4, consider an even element $J_{\fg}$ 
of conformal weight $1$ of the algebra $D^{ch}_{\Omega^{\cdot};\fg}$,  
given by 
$$
J_{\fg}=\phi_{i(-1)}\psi_i=\phi_i\psi_i
\eqno{(6.15.1)}
$$
After a change of frame as in {\it loc. cit.}, we get an element 
$$
J_{\fg'}=\phi'_{i(-1)}\psi'_i=\phi'_i\psi'_i-\gamma(\phi'_i,\psi'_i)
\eqno{(6.15.2)}
$$
where we have again used Theorem 6.4. 
We have 
$$
\phi'_i\psi'_i=g^{-1pi}\phi_p g^{iq}\psi_q=\phi_p\psi_p=J_{\fg}
$$
(see (3.2.3)). 
On the other hand, by (3.4.3)  
$$
\gamma(\phi'_i,\psi'_i)=\gamma(g^{-1pi}\phi_p,g^{iq}\psi_q)=
\delta_{pq}g^{-1pi}\dpar g^{iq}=g^{-1pi}\dpar g^{ip}=
tr(g^{-1}\dpar g)
$$
Thus
$$
J_{\fg'}=J_{\fg}-tr(g^{-1}\dpar g)
\eqno{(6.15.3)}
$$

{\bf 6.16.} Consider an odd element $G_{\fg}$ of conformal 
weight $2$  given by 
$$
G_{\fg}=\psi_{i(-1)}\omega_i
\eqno{(6.16.1)}
$$
In the frame $\fg'$,  
$$
G_{\fg'}=(g^{iq}\psi_q)_{(-1)}(g^{-1si}\omega_s)
$$
Note that $\psi_{q(j)}a=0$ for $j\geq 0$ (everything happens in $D$), hence 
it follows from the commutativity formula (1.3.1) that 
$$
a\psi_q=a_{(-1)}\psi_q=\psi_{q(-1)}a
\eqno{(6.16.2)}
$$
Therefore by "associativity" (1.2.5) 
$$
(a\psi_q)_{(-1)}(b\omega_s)=(\psi_{q(-1)}a)_{(-1)}(b\omega_s)=
\psi_{q(-1)}a_{(-1)}b\omega_s=\psi_{q(-1)}(ab\omega_s)
\eqno{(6.16.3)}
$$
Therefore 
$$
G_{\fg'}=\psi_{q(-1)}g^{iq}g^{-1si}\omega_s=\psi_{s(-1)}\omega_s=G_{\fg}
\eqno{(6.16.6)}
$$

\bigskip\bigskip

{\bf 6.17.} Let us investigate the 
{\it Virasoro element}. Define an even element $L_{\fg}$ of 
conformal weight $2$ by 
$$
L_{\fg}=L_{(b)\fg}+L_{(f)\fg}
\eqno{(6.17.1)}
$$
where
$$
L_{(b)\fg}=\omega_{i(-1)}\tau_i
\eqno{(6.17.2)}
$$
($(b)$ is for "bosonic") and 
$$
L_{(f)\fg}=\rho_{i(-1)}\psi_i
\eqno{(6.17.3)}
$$
($(f)$ is for "fermionic"), cf. [MSV], (2.3a). 

{\bf 6.18.} We have 
$$
(a\omega_s)_{(-1)}(b\tau_p)=
\omega_{s(-1)}\bigl\{ab\tau_p+\tau_p(a)\dpar b+\tau_p(b)\dpar a\bigr\}- 
\dpar\omega_{s(-1)}b\tau_p(a)
\eqno{(6.18.1)}
$$
Indeed, it follows from "associativity" (1.2.5) that 
$$
(a\omega_s)_{(-1)}(b\tau_p)=(\omega_{s(-1)}a)_{(-1)}(b\tau_p)=
\omega_{s(-1)}a_{(-1)}(b\tau_p)+\omega_{s(-2)}a_{(0)}(b\tau_p)+
$$
$$
+a_{(-2)}\omega_{s(0)}(b\tau_p)
$$
Next, 
$$
a_{(-1)}(b\tau_p)=ab\tau_p-\gamma(a,b\tau_p)=
ab\tau_p+\tau_p(a)\dpar b+\tau_p(b)\dpar a
$$
(see [GII] (3.3.1)); 
$$
a_{(0)}(b\tau_p)=-b\tau_{p(0)}a=-b\tau_p(a)
$$
(see [GII] (3.3.2)). Finally 
$$
\omega_{s(0)}(b\tau_p)=-(b\tau_p)_{(0)}\omega_s+\dpar\langle 
b\tau_p,\omega_s\rangle
$$
where
$$
(b\tau_p)_{(0)}\omega_s=(b\tau_p)(\omega_s)=\langle\tau_p,\omega_s\rangle
\dpar b=\delta_{ps}\dpar b
$$
by (1.1.3), since 
$$
\tau_p(\omega_s)=0, 
\eqno{(6.18.2)}
$$
and $\langle b\tau_p,\omega_s\rangle=b\delta_{ps}$. This implies 
$$
\omega_{s(0)}(b\tau_p)=0
\eqno{(6.18.3)}
$$
Formula (6.18.1) follows. 

{\bf 6.19.} We have 
$$
(a\omega_s)_{(-1)}(b\phi_{\alpha}\psi_{\beta})=
\omega_{s(-1)}\bigl\{ab\phi_{\alpha}\psi_{\beta}-\delta_{\alpha\beta}
b\dpar a\bigr\}
\eqno{(6.19.1)}
$$   
Indeed, by commutativity and "associativity" (1.2.5) 
$$    
(a\omega_s)_{(-1)}(b\phi_{\alpha}\psi_{\beta})=
(\omega_{s(-1)}a)_{(-1)}(b\phi_{\alpha}\psi_{\beta})=
\omega_{s(-1)}a_{(-1)}(b\phi_{\alpha}\psi_{\beta})
$$
On the other hand 
$$
a_{(-1)}(b\phi_{\alpha}\psi_{\beta})=
ab\phi_{\alpha}\psi_{\beta}-\gamma(a,b\phi_{\alpha}\psi_{\beta})=
ab\phi_{\alpha}\psi_{\beta}-\delta_{\alpha\beta}b\dpar a, 
$$
see (3.4.2). This implies (6.19.1). 

{\bf 6.20.} We have 
$$
L_{(b)\fg'}=L_{(b)\fg}+\omega_{s(-1)}\bigl\{
\tau_p(g^{-1si})\dpar g^{ip}+g^{-1si}\tau_{\beta}(g^{i\alpha})
\phi_{\beta}\psi_{\alpha}\bigr\}-
$$
$$
-\omega_{s(-2)}g^{ip}\tau_p(g^{-1si})
\eqno{(6.20.1)}
$$
Indeed, due to Theorem 6.4 
$$
L_{(b)\fg'}=(g^{-1si}\omega_s)_{(-1)}\bigl\{
g^{ip}\tau_p+g^{i\alpha\beta}\phi_{\beta}\psi_{\alpha}\bigr\}
$$ 
According to (6.18.1) 
$$
(g^{-1si}\omega_s)_{(-1)}(g^{ip}\tau_p)=
\omega_{s(-1)}\bigr\{g^{-1si}g^{ip}\tau_p+
\tau_p(g^{-1si})\dpar g^{ip}+\tau_p(g^{ip})\dpar g^{-1si}\bigr\}-
$$
$$
-\omega_{s(-2)}g^{ip}\tau_p(g^{-1si})
$$
By (6.18.2)
$$
(g^{-1si}\omega_s)_{(-1)}(g^{i\alpha\beta}\phi_{\beta}\psi_{\alpha})=
(g^{-1si}\omega_s)_{(-1)}(\tau_{\beta}(g^{i\alpha})\phi_{\beta}\psi_{\alpha})=
$$
$$
=\omega_{s(-1)}\bigl\{g^{-1si}\tau_{\beta}(g^{i\alpha})
\phi_{\beta}\psi_{\alpha}-\delta_{\alpha\beta}\tau_{\beta}(g^{i\alpha})
\dpar g^{-1si}\bigr\}
$$
The third term in the first expression cancels out the second term in the 
second one, and we get (6.20.1). 

{\bf 6.21.} We have 
$$
(a\rho_{\mu})_{(-1)}(b\psi_{\nu})=
\rho_{\mu(-1)}ab\psi_{\nu}
\eqno{(6.21.1)}
$$
Indeed, 
$$
(a\rho_{\mu})_{(-1)}(b\psi_{\nu})=
(\rho_{\mu(-1)}a)_{(-1)}(b\psi_{\nu})=
\rho_{\mu(-1)}ab\psi_{\nu}+a_{(-2)}\rho_{\mu(0)}b\psi_{\nu}
$$
and by commutativity  
$$
\rho_{\mu(0)}b\psi_{\nu}=(b\psi_{\nu})(\rho_{\mu})-\dpar\bigl\{
(b\psi_{\nu})_{(1)}\rho_{\mu}\bigr\}=\delta_{\nu\mu}\dpar b-
\delta_{\nu\mu}\dpar b=0,
$$
cf. (6.18.3). This implies (6.21.1).

{\bf 6.22.} We have 
$$
(a\phi_{\gamma}\omega_i)_{(-1)}(b\psi_{\nu})=
\omega_{i(-1)}\bigl\{
ab\phi_{\gamma}\psi_{\nu}-\delta_{\gamma\nu}a\dpar b\bigr\}+
\delta_{\nu\gamma}\omega_{i(-2)}ab
\eqno{(6.22.1)}
$$
Indeed, 
$$
(a\phi_{\gamma}\omega_i)_{(-1)}(b\psi_{\nu})=
(\omega_{i(-1)}a\phi_{\gamma})_{(-1)}(b\psi_{\nu})=
\omega_{i(-1)}(a\phi_{\gamma})_{(-1)}b\psi_{\nu}+
\omega_{i(-2)}(a\phi_{\gamma})_{(0)}b\psi_{\nu}
$$
where
$$
(a\phi_{\gamma})_{(-1)}b\psi_{\nu}=a\phi_{\gamma}b\psi_{\nu}-
\gamma(a\phi_{\gamma},b\psi_{\nu})=
ab\phi_{\gamma}\psi_{\nu}-\delta_{\gamma\nu}a\dpar b, 
$$
cf. (3.4.3), and 
$$
(a\phi_{\gamma})_{(0)}b\psi_{\nu}=
(b\psi_{\nu})(a\phi_{\gamma})=ab\delta_{\nu\gamma}
$$

{\bf 6.23.} We have
$$
L_{(f)\fg'}=L_{(f)\fg}+\omega_{i(-1)}\bigl\{
\tau_i(g^{-1\gamma\alpha})g^{\alpha\nu}\phi_{\gamma}\psi_{\nu}-
\tau_i(g^{-1\gamma\alpha})\dpar g^{\alpha\gamma}\bigr\}+
\omega_{i(-2)}\tau_i(g^{-1\gamma\alpha})g^{\alpha\gamma}
\eqno{(6.23.1)}
$$
Indeed,  
$$
L_{(f)\fg'}=\rho'_{\alpha(-1)}\psi'_{\alpha}=
\bigl\{g^{-1\mu\alpha}\rho_{\mu}+
\tau_i(g^{-1\gamma\alpha})\phi_{\gamma}\omega_i\bigr\}_{(-1)}
(g^{\alpha\nu}\psi_{\nu})
$$
By (6.21.1)
$$
(g^{-1\mu\alpha}\rho_{\mu})_{(-1)}
(g^{\alpha\nu}\psi_{\nu})=\rho_{\mu(-1)}g^{-1\mu\alpha}
g^{\alpha\nu}\psi_{\nu}=L_{(f)\fg}
$$
and by (6.22.1) 
$$
(\tau_i(g^{-1\gamma\alpha})\phi_{\gamma}\omega_i)_{(-1)}
(g^{\alpha\nu}\psi_{\nu})=
\omega_{i(-1)}\bigl\{
\tau_i(g^{-1\gamma\alpha})g^{\alpha\nu}\phi_{\gamma}\psi_{\nu}-
\tau_i(g^{-1\gamma\alpha})\dpar g^{\alpha\gamma}\bigr\}+
$$
$$
+\omega_{i(-2)}\tau_i(g^{-1\gamma\alpha})g^{\alpha\gamma}
$$

{\bf 6.24.} Comparing (6.20.1) and (6.23.1) we see easily that 
$$
L_{\fg'}=L_{(b)\fg'}+L_{(f)\fg'}=
L_{(b)\fg}+L_{(f)\fg}=L_{\fg}
\eqno{(6.24.1)}
$$
Let us collect our computations of transformation rules. 

{\bf 6.25. Theorem.} {\it Let $\fg, \fg'$ be two natural frames of 
$(A,\Omega)$.  
Consider $4$ elements of the 
vertex superalgebra $D^{ch}_{\Omega^{\cdot};\fg}$ given by 
$$
Q_{\fg}=\phi_{i(-1)}\tau_i
\eqno{(6.25.1)}
$$
(an odd element of conformal weight $1$)
$$
J_{\fg}=\phi_{i(-1)}\psi_i
\eqno{(6.25.2)}
$$
(an even element of conformal weight $1$)
$$
G_{\fg}=\psi_{i(-1)}\omega_i
\eqno{(6.25.3)}
$$
(an odd element of conformal weight $2$) and
$$
L_{\fg}=\omega_{i(-1)}\tau_i+\rho_{i(-1)}\psi_i
\eqno{(6.25.4)}
$$
(an even element of conformal weight $2$). 

After the canonical identification $D^{ch}_{\Omega^{\cdot};\fg'}=
D^{ch}_{\Omega^{\cdot};\fg}$ these elements are transformed as follows
$$
Q_{\fg'}=Q_{\fg}+\dpar\bigl\{
tr\bigl(g^{-1}\tau_r(g)\bigr)\phi_r\bigr\}
\eqno{(6.25.5)}
$$
$$
J_{\fg'}=J_{\fg}-tr(g^{-1}\dpar g)
\eqno{(6.25.6)}
$$
$$
G_{\fg'}=G_{\fg}
\eqno{(6.25.7)}
$$
and
$$
L_{\fg'}=L_{\fg}
\eqno{(6.25.8)}
$$}

This is a version of [MSV], Theorem 4.2.

\bigskip\bigskip 

\centerline{\bf \S 7. Poincar\'e-Birkhoff-Witt}

\bigskip\bigskip

{\bf 7.1} Let $X$ be a smooth variety and $\CD^{ch}_{\Omega^{\cdot}_X}$ 
be the sheaf discussed in the previous section, cf. Theorem 6.6. It is 
a sheaf of $\BZ_{\geq 0}$-graded vertex algebras, so 
$$
\CD^{ch}_{\Omega^{\cdot}_X}=\oplus_{n\in \BZ_{\geq 0}}\ 
\CD^{ch}_{\Omega^{\cdot}_X; n}
\eqno{(7.1.1)}
$$
where $\CD^{ch}_{\Omega^{\cdot}_X; n}$ denotes the component of conformal 
weight $n$. 

According to Theorem 6.25 (see (6.25.8)) we have a canonical 
global section $L$ of $\CD^{ch}_{\Omega^{\cdot}_X; 2}$. Let $L(z)=
\sum\ L_nz^{-n-2}$ be the corresponding field.  
 
{\bf 7.1.1. Claim.} {\it A local section 
$\alpha\in \CD^{ch}_{\Omega^{\cdot}_X}$ belongs to 
$\CD^{ch}_{\Omega^{\cdot}_X; n}$ if and only if $L_0(\alpha)=n\alpha$.} 

In a uniform notation $L(z)=\sum\ L_{(n)}z^{-n-1}$ we have $L_0=L_{(1)}$. 
We shall check 7.1.1 simultaneously with 

{\bf 7.1.2. Claim.} {\it The operator $L_{-1}=L_{(0)}$ coincides with 
the canonical derivation $\dpar$ of the vertex algebra 
$\CD^{ch}_{\Omega^{\cdot}_X}$.} 

Both statements are local, so we may assume we are in the local situation 
6.17. Note that the operator $L_{(0)}$ is a derivation with respect 
to the operation $_{(-1)}$: 
$$
L_{(0)}(x_{(-1)}y)=(L_{(0)}x)_{(-1)}y+x_{(-1)}L_{(0)}y
\eqno{(7.1.2)}
$$
Therefore it suffices to check 7.1.2 on the generators $a, \tau_i, \omega_i, 
\psi_i, \rho_i$ of the vertex algebroid $\CA_{\Omega^{\cdot};\fg}$, 
which is done by a simple explicit computation. 

It follows from the associativity formula (1.2.4) that 
$$
L_{(1)}y_{(-1)}z=(L_{(1)}y)_{(-1)}z+y_{(-1)}L_{(1)}z+L_{(0)}y_{(0)}z 
-y_{(0)}L_{(0)}z=
$$
$$
=(L_{(1)}y)_{(-1)}z+y_{(-1)}L_{(1)}z
\eqno{(7.1.3)}
$$
since
$$       
L_{(0)}y_{(0)}z - y_{(0)}L_{(0)}z=\dpar(y_{(0)}z)-y_{(0)}\dpar z=0
$$
In other words, $L_{(1)}$ is a derivation of $_{(-1)}$. Therefore it 
suffices to check 7.1.1 on the generators of $\CA_{\Omega^{\cdot};\fg}$ 
as above, which is straightforward.

{\bf 7.2.} In the local situation 6.1, consider the local algebra 
$D^{ch}_{\Omega^{\cdot};\fg}=U\CA_{\Omega^{\cdot};\fg}$. Let us introduce 
a $\BZ$-grading 
$$
D^{ch}_{\Omega^{\cdot};\fg}=\oplus_{p\in Z}\ 
D^{ch; p}_{\Omega^{\cdot};\fg}
\eqno{(7.2.1)}
$$
to be called {\it fermionic charge}. For an element $x\in 
D^{ch}_{\Omega^{\cdot};\fg}$ let us denote its fermionic charge (to be 
defined) by $F(x)\in\BZ$. It is defined uniquely by the following 
conditions: 

(a) $F(a)=F(\tau_i)=F(\omega_i)=0;\ F(\phi_i)=F(\rho_i)=-F(\psi_i)=1$; 

(b) $F(x_{(-1)}y)=F(x)+F(y)$.  

Due to the transformation formulas (3.2.1) - (3.2.5) this grading is 
obviously preserved under a change of frames. Therefore in the situation 
of 7.1 the sheaf $\CD^{ch}_{\Omega^{\cdot}_X}$ gets a canonical 
$\BZ$-grading 
$$
\CD^{ch}_{\Omega^{\cdot}_X}=\oplus_{p\in\BZ}\ 
\CD^{ch; p}_{\Omega^{\cdot}_X}
\eqno{(7.2.2)}
$$
Note that parity 
is equal to fermionic charge modulo $2$. 

Here is another way to define the grading (7.2.2). First 
notice a simple 

{\bf 7.2.1. Lemma.} {\it Let $\CA=(A,T,\Omega,\dpar,\ldots)$ be a 
vertex (super)algebroid. For every invertible element $a\in A$ the operator 
$(a^{-1}\dpar a)_{(0)}$ acting on $U\CA$ is trivial.} 

{\it Proof.} Obviously this operator is trivial on $A=U\CA_0$. Let 
$x\in U\CA_1$ and $\tau\in T$ be its image under the canonical projection 
$U\CA_1\lra T$. We have 
$$
(a^{-1}\dpar a)_{(0)}x=-x_{(0)}a^{-1}\dpar a+\dpar(x_{(1)}a^{-1}\dpar a)=
-\tau(a^{-1}\dpar a)+\dpar\langle\tau,a^{-1}\dpar a\rangle=
$$
$$
=a^{-2}\tau(a)\dpar a-a^{-1}\dpar \tau(a)+\dpar(a^{-1}\tau(a))=0,
$$   
so $(a^{-1}\dpar a)_{(0)}$ is trivial on $U\CA_1$. Therefore it is 
trivial on the whole algebra $U\CA$ since $(?)_{(0)}$ is a derivation 
of the operation $_{(-1)}$. $\btu$

Applying this lemma to $a=det(g)$ in the formula (6.25.6) we see that  
the component $J_{0;\fg}$ of the field $J_{\fg}(z)=
\sum\ J_{n;\fg}z^{-n-1}$ is preserved under the change of frames. 
Consequently it gives rise to a well defined endomorphism $J_0$ 
of the sheaf $\CD^{ch}_{\Omega^{\cdot}_X}$.   

{\bf 7.2.2. Claim.} {\it A local section 
$\alpha\in \CD^{ch}_{\Omega^{\cdot}_X}$ belongs to 
$\CD^{ch; p}_{\Omega^{\cdot}_X}$ if and only if $J_0(\alpha)=p\alpha$.} 

Indeed, the function $F(\alpha)$ defined by $J_0(\alpha)=F(\alpha)\alpha$ 
obviously satisfies the condtions (a) and (b) above.  

{\bf 7.3.} The two gradings (7.1.1) and (7.2.2) are compatible: if we denote 
$$
\CD^{ch; p}_{\Omega^{\cdot}_X; n}:=
\CD^{ch}_{\Omega^{\cdot}_X; n}\cap \CD^{ch; p}_{\Omega^{\cdot}_X}
\eqno{(7.3.1)}
$$
then 
$$
\CD^{ch}_{\Omega^{\cdot}_X}=\oplus_{n\in\BZ_{\geq 0};\ p\in \BZ}\ 
\CD^{ch; p}_{\Omega^{\cdot}_X; n}
\eqno{(7.3.2)}
$$
For a fixed $n$, only a finite number of sheaves 
$\CD^{ch; p}_{\Omega^{\cdot}_X; n}$ are nonzero. 

If the ground ring $k$ is a field of characteristic $0$ then the sheaves
$\CD^{ch; p}_{\Omega^{\cdot}_X; n}$ and 
$\CD^{ch; n-p}_{\Omega^{\cdot}_X; n}$ are in a certain sense dual to each 
other,   
see [MS]. 

{\bf 7.4.} Starting from this point we assume that $k\supset\BQ$. 
According to a (superversion of) the PBW theorem, [GII], Theorem 9.18, 
the sheaf $\CD^{ch}_{\Omega^{\cdot}_X}$ admits a canonical filtration 
such that the associated graded sheaf is canonically isomorphic to 
$$
gr(\CD^{ch}_{\Omega^{\cdot}_X})=Sym_{\Omega^{\cdot}_X}\bigl\{ 
\bigl(\oplus_{n\geq 1}\ \Theta_{\Omega^{\cdot}_X (n)}\bigr)\oplus
\bigl(\oplus_{n\geq 1}\ \Omega^1_{\Omega^{\cdot}_X (n)}\bigr)\bigr\}
\eqno{(7.4.1)}
$$
Here $\Theta_{\Omega^{\cdot}_X}$ (resp. $\Omega^1_{\Omega^{\cdot}_X}$) 
denotes the tangent (resp. the cotangent) bundle of the supervariety 
$(X,\Omega^{\cdot}_X)$, and $(?)_{(n)}$ means that this bundle is put 
into the conformal weight $n$. 

The endomorphisms $L_0$ and $J_0$ respect the canonical filtration; 
hence we get a canonical finite filtration $F_{\cdot}$ on each homogeneous 
component $\CD^{ch; p}_{\Omega^{\cdot}_X; n}$. The graded quotients 
$F_i\CD^{ch; p}_{\Omega^{\cdot}_X; n}/F_{i+1}\CD^{ch; p}_{\Omega^{\cdot}_X; n}$ 
are locally free $\CO_X$-modules of finite rank (we shall see this 
in the course of computations below). This allows us to introduce the 
elements of the Grothendieck group $K(X)$ of vector bundles 
$$
[\CD^{ch; p}_{\Omega^{\cdot}_X; n}]:=
\sum_i\ [F_i\CD^{ch; p}_{\Omega^{\cdot}_X; n}/
F_{i+1}\CD^{ch; p}_{\Omega^{\cdot}_X; n}]\in K(X)
\eqno{(7.4.2)}
$$
Here $[E]$ in the right hand side denotes the class of a vector bundle 
$E$ in $K(X)$. 
Consider the generating function 
$$
cl(\CD^{ch}_{\Omega^{\cdot}_X})(y,q):=
\sum_{p, n}\ [\CD^{ch; p}_{\Omega^{\cdot}_X; n}] y^p q^n\in 
K(X)[y,y^{-1}][[q]]
\eqno{(7.4.3)}
$$         

{\bf 7.5.} For a vector bundle $E$ over $X$ and an indeterminate $x$ we 
introduce the notations 
$$
[S_x E]=\sum_{i=0}^{\infty}\ [Sym^i_{\CO_X} E]\in K(X)[[x]]
\eqno{(7.5.1)}
$$
and 
$$
[\Lambda_x E]=\sum_{i=0}^{\infty}\ 
[\Lambda^i_{\CO_X} E]\in K(X)[x]
\eqno{(7.5.2)}
$$
The following fact was noticed in [BL] (cf. also [W]).  

{\bf 7.6. Theorem} (L. Borisov - A. Libgober) {\it We have 
$$
cl(\CD^{ch}_{\Omega^{\cdot}_X})(y,q)=[\Lambda_y\Omega^1_X]\cdot 
\prod_{n=1}^{\infty}\bigl\{ 
[S_{q^n}\Theta_X]\cdot [S_{q^n}\Omega^1_X]\cdot
[\Lambda_{y^{-1}q^n}\Theta_X]\cdot [\Lambda_{yq^n}\Omega^1_X]\bigr\}  
\eqno{(7.6.1)}
$$} 

{\bf 7.7.} {\it Proof.} Let us understand the bundles 
$\Theta_{\Omega^{\cdot}_X}$ and $\Omega^1_{\Omega^{\cdot}_X}$ 
a little bit more attentively.   

Let us consider the local situation 3.1, 
with $E=\Omega$, so that $\Lambda E=\Omega^{\cdot}$.  
All our frames $\fg$ will be natural. Let 
$T_{\psi}\subset T_{\Omega^{\cdot}}$ be the $A$-submodule with the base 
$\{\psi_i\}$. The coordinate change formula (3.2.3) shows that 
it is a well defined $A$-submodule of $T$ independent on the 
choice of a frame, canonically isomorphic to $T$. 

We set 
$$
T_{\psi\Omega^{\cdot}}=\Omega^{\cdot}\otimes_A T_{\psi}\subset 
T_{\Omega^{\cdot}}
$$
We denote by $T_{\tau\Omega^{\cdot}}$ the quotient $\Omega^{\cdot}$-module 
$T_{\Omega^{\cdot}}/T_{\psi\Omega^{\cdot}}$. Let $T_{\tau}\subset 
T_{\tau\Omega^{\cdot}}$ be the $A$-submodule generated by all $\tau_i$. 
The formula (3.2.1) shows that $T_{\tau}$ is a well defined $A$-module 
canonically isomorphic to $T$, and we have 
$$
T_{\tau\Omega^{\cdot}}=\Omega^{\cdot}\otimes_A T_{\tau}
$$
Returning to our variety $X$, we see that we get two vector bundles  
$\Theta_{\psi}$ and $\Theta_{\tau}$ both isomorphic to $\Theta_X$ and 
a canonical short exact sequence 
$$
0\lra \Theta_{\psi\Omega^{\cdot}_X}\lra \Theta_{\Omega^{\cdot}_X}\lra 
\Theta_{\tau\Omega^{\cdot}_X}\lra 0
\eqno{(7.7.1)}
$$
with
$$
\Theta_{\psi\Omega^{\cdot}_X}=\Omega^{\cdot}_X\otimes_{\CO_X} 
\Theta_{\psi};\ \    
\Theta_{\tau\Omega^{\cdot}_X}=\Omega^{\cdot}_X\otimes_{\CO_X} 
\Theta_{\tau} 
\eqno{(7.7.2)}
$$
Note that $\Theta_{\psi}$ has fermionic charge $-1$ and 
$\Theta_{\tau}$ has fermionic charge $0$. 

Dually, we have two vector bundles  
$\Omega^1_{\rho}$ and $\Omega^1_{\omega}$ both isomorphic to $\Omega^1_X$ and 
a canonical short exact sequence 
$$
0\lra \Omega^1_{\omega\Omega^{\cdot}_X}\lra \Omega^1_{\Omega^{\cdot}_X}\lra 
\Omega^1_{\rho\Omega^{\cdot}_X}\lra 0
\eqno{(7.7.3)}
$$
with
$$
\Omega^1_{\rho\Omega^{\cdot}_X}=\Omega^{\cdot}_X\otimes_{\CO_X} 
\Omega^1_{\rho};\ \    
\Omega^1_{\omega\Omega^{\cdot}_X}=\Omega^{\cdot}_X\otimes_{\CO_X} 
\Omega^1_{\omega} 
\eqno{(7.7.4)}
$$
The bundles $\Omega^1_{\rho}$, $\Omega^1_{\omega}$ have fermionic 
charges $1$, $0$ respectively. 

Note that if $E$ is a vector bundle then 
$$
Sym_{\Omega^{\cdot}_X}(\Omega^{\cdot}_X\otimes_{\CO_X}E)=
\Omega^{\cdot}_X\otimes_{\CO_X} Sym_{\CO_X} E
\eqno{(7.7.5)}
$$ 
Returning to PBW formula (7.4.1) we see that these remarks imply       
(7.6.1). $\btu$  

{\bf 7.8.} Starting from this point let us assume that $k=\BC$. 
Consider the formal power series 
$$
\theta(y,q)=i^{-1}(y^{1/2}-y^{-1/2})q^{1/8}\prod_{n=1}^{\infty}\bigl\{  
(1-q^n)(1-yq^n)(1-y^{-1}q^n)\bigr\}
\eqno{(7.8.1)}
$$
It is nothing but the theta function $\theta_1(h,z)$ as defined in 
[HC], II, 2, \S 10, formula (3), p. 204, with $q=h^2$ and $y=z^2$. 

If $f\in GL(V)$ is an automorphism of a $d$-dimendional vector space $V$  
with eigenvalues $\lambda_1,\ldots,\lambda_d$, we shall denote by 
$\theta_f(y,q)$ the power series 
$$
\theta_f(y,q)=\frac{\prod_{i=1}^d\ \theta(\lambda_iy,q)}
{\prod_{i=1}^d\ \theta(\lambda_i,q)}
\eqno{(7.8.2)}
$$
Let $X$ be a proper smooth $d$-dimensional algebraic variety; let 
$g:\ X\lra X$ be a 
{\it simple} automorphism, which means by definition that 
the graph $\Gamma_g\subset X\times X$ is transversal to the diagonal. 
This implies that the set $X^g$ if its fixed points is finite. 

For each $x\in X^g$ denote by $g_x$ the induced endomorphism of the 
cotangent space $\Omega^1_{X;x}$. All eigenvalues of $g_x$ are 
distinct from $1$. 

{\bf 7.9. Theorem.} {\it Consider the power series 
$$
T_{X;g}(y,q):=y^{-d/2}\sum_{a,b,n}\ (-1)^{a+b}Tr(g; H^a(X;
\CD^{ch;b}_{\Omega^{\cdot}_X; n}))y^bq^n
\eqno{(7.9.1)}
$$
We have
$$
T_{X;g}(y,q)=\sum_{x\in X^g}\ \theta_{g_x}(y,q)
\eqno{(7.9.2)}
$$}

{\bf 7.10.} {\it Proof.} Recall that according to the Atiyah-Bott holomorphic 
Lefschetz fixed point formula, if $E$ is a $g$-equivariant vector bundle 
over $X$ then 
$$
\sum_i\ (-1)^i Tr(g; H^i(X;E))=\sum_{x\in X^g}\ 
\frac{Tr(g; E_x)}{det(1-g_x)}
\eqno{(7.10.1)}
$$
see [AB], Theorem 4.12. Note that if $(V, f)$ are as in 7.8 then 
$$
Tr(g; Sym_x(V))=\prod_{i=1}^d\ (1-\lambda_ix)^{-1}=
Tr(g;\Lambda_{-x}(V))^{-1}
\eqno{(7.10.2)}
$$
The proof of 7.6 shows that each sheaf 
$\CD^{ch; p}_{\Omega^{\cdot}_X; n}$ carries a canonical filtration 
whose quotients are vector bundles and the associated graded sheaf 
is given by the formula (7.6.1). Therefore we may apply the 
Lefschetz formula (7.10.1). 

Note that since in the expression (7.9.1) the fermionic charge $a+b$ 
is taken into account, we should apply the Lefschetz formula to the 
element $cl(\CD^{ch}_{\Omega^{\cdot}_X})(-y,q)$. Due to (7.10.2) each 
fixed point $x$ gives a contribution 
$$
y^{-d/2}\prod_{i=1}^d\biggl\{
\frac{1-\lambda_iy}{1-\lambda_i} 
\prod_{n=1}^{\infty}
\frac{(1-\lambda_iyq^n)(1-\lambda_i^{-1}y^{-1}q^n)}{(1-\lambda_iq^n)
(1-\lambda_i^{-1}q^n)}
\biggr\}=\theta_{g_x}(y,q)
$$
where $\lambda_i$ are the eigenvalues of $g_x$. This implies the theorem. 
$\btu$. 

The reader may wish to compare (7.9.2) with the explicit formulas for 
the trace 
of certain automorphisms of the Frenkel-Lepowsky-Meurman Monster 
vertex algebra, cf. [FLM].

\bigskip\bigskip

\centerline{\bf References}

\bigskip\bigskip

[AB] M.F.~Atiyah, R.~Bott, A Lefschetz fixed point formula for 
elliptic complexes: II. Applications, {\it Ann. Math.}, {\bf 88}, No. 3 
(1968), 451-491. 

[BL] L.~Borisov, A.~Libgober, Elliptic genera and applications to 
mirror symmetry, math.AG/9904126. 

[FLM] I.~Frenkel, J.~Lepowsky, A.~Meurman, Vertex operator algebras 
and the Monster, {\it Pure and Applied Mathematics}, {\bf 134}, 
Academic Press, Boston, 1988.  

[GI] V.~Gorbounov, F.~Malikov, V.~Schechtman, Gerbes of chiral differential 
operators, math.AG/9906117; {\it Math. Research Letters}, {\bf 7}, 
1-12 (2000). 

[GII] V.~Gorbounov, F.~Malikov, V.~Schechtman, Gerbes of chiral 
differential operators. II, math.AG/0003170. 

[HC] A.~Hurwitz, R.~Courant, Funktionentheorie, Vierte Auflage, 
Springer-Verlag, Berlin-G\"ottingen-Heidelberg-New York, 1964.  

[K] V.~Kac, Vertex algebras for beginners, 
Second Edition, University Lecture Series, {\bf 10}, 
American Mathematical Society, Providence, Rhode Island, 1998. 

   
[MSV] F.~Malikov, V.~Schechtman, A.~Vaintrob, Chiral de Rham complex,\ 
{\it Comm. Math. Phys.}, {\bf 204} (1999), 439-473. 

[MS] F.~Malikov, V.~Schechtman, Chiral Poincar\'e duality, 
{\it Math. Research Letters}, {\bf 6} (1999), 533-546. 

[W] E.~Witten, Elliptic genera and quantum field theory, 
{\it Comm. Math. Phys.} {\bf 109} (1987), 525-536.

\bigskip

\bigskip

V.G.: Department of Mathematics, University of Kentucky, 
Lexington, KY 40506, USA;\ vgorb\@ms.uky.edu

F.M.: Department of Mathematics, University of Southern California, 
Los Angeles, CA 90089, USA;\ fmalikov\@mathj.usc.edu   

V.S.: IHES, 35 Route de Chartres, 91440 Bures-sur-Yvette, France;\ 
vadik\@ihes.fr

\enddocument